\documentclass[journal]{IEEEtran}
\ifCLASSINFOpdf
\else
\fi
\usepackage{amsmath}

\usepackage{algorithm}
\usepackage{algpseudocode}
\usepackage{units}
\usepackage{prettyref} 
\newrefformat{sec}{\mbox{Section \ref{#1}}}
\newrefformat{tab}{\mbox{Tab. \ref{#1}}}
\newrefformat{fig}{\mbox{Fig. \ref{#1}}}
\newrefformat{algo}{\mbox{Algorithm \ref{#1}}}
\usepackage{epsfig} 
\usepackage{booktabs}
\usepackage{graphics} 
\usepackage{xcolor}
\usepackage{bbm}
\usepackage[super,negative]{nth}
\usepackage{algpseudocode}
\usepackage{enumitem}
\usepackage{amsthm}
\usepackage{thmtools}

\usepackage{nicefrac}

\usepackage{url}
\usepackage{amssymb}
\usepackage{soul}
\usepackage{graphicx}

\newcommand{\T}{\scriptscriptstyle\top}       
\newcommand{\mathmin}{\operatorname*{min}}
\newcommand{\mathst}{\text{s.t.}}

\definecolor{myOrange}{RGB}{255,165,0}
\definecolor{myGreen}{rgb}{0.44314,0.74902,0.43137}%
\definecolor{mycolor2}{rgb}{0.91373,0.28235,0.28627}%

\newcommand{\chgAlex}[1]{#1}
\newcommand{\chgRevI}[1]{#1}
\newcommand{\chgRevII}[1]{#1}

\newcommand{\chgRevIV}[1]{#1}

\declaretheoremstyle[headfont=\bfseries]{normalhead}
\declaretheorem[style=normalhead]{assumption}
\declaretheorem[style=normalhead]{remark}
\declaretheorem[name={Control Problem},style=normalhead]{problem}

\usepackage{tikz}
\usetikzlibrary{shapes,arrows}
\usepackage{pgfplots} 
\usepackage{pgfgantt}
\usepackage{pdflscape}
\pgfplotsset{compat=1.5}
\pgfplotsset{plot coordinates/math parser=false}
\newlength\fwidth
\usepackage{grffile}
\usetikzlibrary{plotmarks}
\usetikzlibrary{arrows.meta}
\usepgfplotslibrary{patchplots}

\begin{document}
%
\title{Fully Distributed Model Predictive Control of Connected Automated Vehicles in Intersections: Theory and Vehicle Experiments} 
%
%
%

\author{Alexander~Katriniok,~\IEEEmembership{Senior Member,~IEEE}, Benedikt~Rosarius~and~Petri~M\"{a}h\"{o}nen,~\IEEEmembership{Senior Member,~IEEE}
	\thanks{A. Katriniok is with the Ford Research \& Innovation Center, 52072 Aachen, Germany, {\tt\small de.alexander.katriniok@ieee.org}.}
	\thanks{B. Rosarius was with the Ford Research \& Innovation Center, 52072 Aachen, Germany, {\tt\small benedikt.rosarius@rwth-aachen.de}.}
	\thanks{P. M\"{a}h\"{o}nen is with the Institute for Networked Systems (INETS), Department of Electrical Engineering, RWTH Aachen University, 52072 Aachen, Germany, {\tt\small pma@inets.rwth-aachen.de}.}}

%
%

\markboth{IEEE Transactions on Intelligent Transportation Systems,~ Katriniok \MakeLowercase{\textit{et al.}}}%
{Katriniok \MakeLowercase{\textit{et al.}}: Bare Demo of IEEEtran.cls for IEEE Journals}
%



\newcommand\copyrighttext{%
	\footnotesize \textcopyright \,2021 IEEE. Personal use of this material is permitted. Permission from IEEE must be obtained for all other uses, in any current or future media, including reprinting/republishing this material for advertising or promotional purposes, creating new collective works, for resale or redistribution to servers or lists, or reuse of any copyrighted component of this work in other works.}
\newcommand\copyrightnotice{%
	\begin{tikzpicture}[remember picture,overlay]
	\node[anchor=south,yshift=4pt] at (current page.south) {\parbox{\dimexpr\textwidth-\fboxsep-\fboxrule\relax}{\copyrighttext}};
	\end{tikzpicture}%
}

\maketitle

\copyrightnotice

\vspace*{-2.5mm}
\begin{abstract}
We propose a fully distributed control system architecture, amenable to in-vehicle implementation, that aims to safely coordinate connected and automated vehicles (CAVs) \chgRevIV{at} road intersections. For control purposes, we build upon a fully distributed model predictive control approach, in which the agents solve a nonconvex optimal control problem (OCP) locally and synchronously, and exchange their optimized trajectories via vehicle-to-vehicle (V2V) communication. To accommodate a fast solution of the nonconvex OCPs, we apply the penalty convex-concave procedure which \chgAlex{solves} a convexified version of the original OCP. For experimental evaluation, we complement the predictive controller with a localization layer, being in charge of self-localization, \chgAlex{and an estimator, which determines joint collision points with other agents}. Experimental tests reveal the efficacy of the proposed control system architecture.
\end{abstract}

\vspace*{-1mm}


%
\IEEEpeerreviewmaketitle

\vspace*{-2mm}
\section{Introduction}
\label{sec:introduction}
\IEEEPARstart{T}{he} automation of road vehicles utilizing vehicle-to-everything (V2X) communication is an emerging field and will support many advancements in intelligent transportation systems \cite{Wymeersch2015}. Connected and automated vehicles (CAVs) are equipped with a communication device for mutual data exchange with other vehicles, the infrastructure or even vulnerable road users. This ability to communicate may complement conventional on-board sensors like radars, cameras or LiDARs such as to extend \chgAlex{their} sensing capabilities in terms of range and the detection of occluded objects. With a sufficient penetration in the market, CAVs can even operate more proactively by negotiating control actions instead of reacting on instantaneous measurements or predicted (but still uncertain) motion trajectories of surrounding vehicles. Potential use cases may involve, amongst others, collaborative lane change maneuvers or \chgAlex{the automation of road intersections} \cite{Khayatian2020a}. In this article, we particularly aim to address \chgAlex{the latter problem, that is, to safely coordinate} CAVs \chgRevIV{at road intersections with no traffic signs or lights}.

\subsection{\chgRevI{Related Work}}
\label{sec:introduction_relatedWork}

For the problem at hand, there is a rich body of literature. Very recent and comprehensive surveys in that space can be found in \cite{Khayatian2020a,Namazi2019a}. From an architectural viewpoint, the respective control schemes can be categorized into centralized, distributed, decentralized and hybrid approaches. Centralized control regimes \chgRevI{\cite{Quinlan2010,Kamal2015,Murgovski2015,Mueller2016a}} require the vehicles (also referred to as \textit{agents}) to communicate with a central node, which then grants exclusive access to the intersection or optimizes the agents' trajectories through the intersection. In distributed schemes \chgRevI{\cite{Makarem2013,Campos2014,Katriniok2017a,Malikopoulos2021a}}, the agents communicate with each other and solve their part of the control problem locally without the involvement of any central node. Decentralized approaches \chgRevI{\cite{Medina2015,Wu2019a,Tian2020a,Schildbach2016a}} differ from distributed concepts in a sense that they do not even involve any communication. Finally, hybrid approaches \cite{Kim2014,Hult2016,Gregoire2016a} are a combination of the aforementioned architectures, that is, these combine, e.g., a centralized regime, assigning a passing order to the agents, with distributed or decentralized controllers being in charge of determining appropriate control actions. 
Contemplating the applied methodology, the intersection coordination problem has amongst others been addressed through hybrid system theory 
\cite{Gregoire2016a,Hafner2013}, \chgRevI{responsibility-sensitive safety (RSS) rules \cite{Khayatian2021a},} resource reservation protocols \cite{Dresner2008,Kowshik2011}, scheduling-based approaches \cite{Ahn2014,Colombo2015}, game theory \cite{Tian2020a,Wei2018a}, \chgRevI{virtual platoons \cite{Medina2015,Englund2016a},} reinforcement learning \cite{Wu2019a,Lamouik2017a} or optimization-based \mbox{control \cite{Kamal2015,Katriniok2017a,Kim2014,Hult2019a,Liu2018a}}.

Narrowing our focus to optimization-based control, \cite{Kamal2015} introduces a centralized model predictive control (MPC) scheme which minimizes the total quantified risk of collision between agents. 
Another centralized MPC scheme is proposed in \cite{Murgovski2015}, where the optimal solution is obtained by solving optimal control subproblems for all combinations of agent crossing sequences. The subproblems are convex and are formulated in the spacial instead of the time domain.
A hybrid approach with a centralized coordination layer, which prescribes the intersection crossing order, and a distributed MPC-based motion planner is outlined in \cite{Kim2014}. A similar idea is pursued in \cite{Kneissl2018} where a central node is in charge of time slot allocation, while agents are controlled in a decentralized fashion. 
The authors in \cite{Campos2014} propose a distributed MPC scheme in which the agents decide sequentially (for a given and fixed decision order) whether to pass the intersection before or after the agents with higher decision order by solving two convex quadratic programming (QP) problems. Another sequential approach within a distributed MPC framework is presented in \cite{Qian2015}. The agents solve their optimal control problem (OCP) sequentially for an \textit{a priori} fixed intersection crossing order. \cite{Liu2018a} suggests to decompose the control problem into a distributed decision maker, which determines intersection entry and exit times, and a distributed motion planner, which optimizes every agent's speed profile such as to meet the respective entry and exit times. Moreover, \cite{Molinari2018} presents a decentralized consensus-based control strategy which determines the intersection crossing order as part of a high-level consensus algorithm and solves a distributed OCP on a lower level to determine vehicle controls. Instead of using time slots, collision avoidance is ensured by imposing a lower bound on agent distances. The hierarchical distributed control scheme in \cite{Malikopoulos2021a} optimizes its optimal trajectory and lane to pass on a higher control level while the optimal vehicle acceleration is determined on a lower control level. A semi-distributed control regime is proposed in \cite{Hult2016}. It utilizes a central coordinator to solve a nonlinear high-level time slot allocation problem for a fixed intersection crossing order, while control actions are determined locally by the agents as part of a nested low-level OCP. Essentially, every agent solves a QP and two linear programming (LP) problems, and transmits their solution to the coordinator which solves a nonlinear programming (NLP) problem. This work is extended towards rear-end collision avoidance in \cite{Shi2018}. 
In \cite{Hult2019a}, the authors present in-vehicle experiments for the use case of straight crossing agents. 
For localization, a centimeter-precision real-time kinematic (RTK) system is utilized.
 
As part of our own research, we have outlined a fully distributed MPC scheme in \cite{Katriniok2017a} where every agent solves a nonconvex quadratically constrained QP (QCQP) through semidefinite relaxation (SDR) with randomization. To decompose the OCP, we introduce \textit{a priori} fixed agent priorities which release the higher priority agent from imposing collision avoidance constraints.
In \cite{Katriniok2019a}, we investigate a reformulation of \cite{Katriniok2017a} which is solved by exploiting a first order optimization method. 
Our work in \cite{Molinari2020a} extends \cite{Katriniok2017a,Molinari2018} towards a hierarchical distributed control architecture which accommodates time-varying agent priorities.
For human driven vehicles, we have proposed a stochastic distributed control regime which issues speed recommendations to the driver \cite{Katriniok2019b}.

\subsection{Main Contribution and Outline}
\label{sec:introduction_contribution}
In this article, we propose a fully distributed control system architecture, which is amenable to in-vehicle implementation and able to safely coordinate CAVs \chgRevIV{at road intersections with no traffic signs or lights}. For control purposes, we utilize a fully distributed MPC scheme, in which the agents solve their respective OCPs synchronously and fast \chgRevIV{to} meet \chgRevIV{real-time} requirements. We consider MPC to be an appealing methodology to approach the control problem at hand as it allows us to explicitly accommodate constraints and to exploit anticipated trajectories of other agents. Compared to a central node, the distributed scheme is more resilient against a single point of failure and scales much better with the number of agents. \chgRevIV{For information exchange among the agents}, 
we rely on vehicle-to-vehicle (V2V) communication, more particularly on Dedicated Short Range Communication (DSRC). 

\chgRevII{We have built upon our previous work \cite{Katriniok2017a} which provided an initial proof of concept under simplifying assumptions in a simulation environment. That said, the real-time implementation in a test vehicle goes far beyond running our algorithm in \cite{Katriniok2017a} on an embedded hardware. Conversely, we had to change and improve our control concept and come up with additional algorithms which are essential for in-vehicle implementation:
\begin{itemize}
	\item \chgRevII{In \cite{Katriniok2017a}, we solve the OCP via SDR with randomization, which is computationally prohibitive on an embedded hardware. To this end, we propose a tailored version of the penalty convex-concave-procedure (CCP) \cite{Lipp2016} which is computationally efficient for embedded implementation and allows the utilization of mature QP solvers.} 
	\item We introduce a self-localization algorithm which accurately and smoothly estimates agent positions by integrating inertial measurements with GNSS measurements. That way, we accommodate the (simplifying) assumption in \cite{Katriniok2017a} that the entire state vector is measurable. 
	\item Another (simplifying) assumption in \cite{Katriniok2017a} is the \textit{a priori} knowledge of joint collision points. For the experimental setup, we design an estimator within the localization layer to determine these collision points online.
	\item The real-time control system is finally implemented on a dSPACE MicroAutoBox II, integrated in our test vehicles and evaluated in experimental tests on the proving ground. Along these lines, we ensure the synchronous execution of local MPC controllers and propose a proprietary V2V communication protocol for information exchange between these controllers.	  
	\item We investigate the feasibility of using a low-cost Global Navigation Satellite System (GNSS) instead of a centimeter-precision RTK as in \cite{Hafner2013,Hult2019a}.
\end{itemize}
In literature, several experiments have already addressed the problem of autonomous intersection crossing, e.g., as part of the DARPA Urban Challenge (DUC) 2017 \cite{Buehler2009a}. The AVs in the DUC, though, had to solve the problem without inter-vehicle communication but through anticipation of the other road users' behavior. Experiments that involve V2V communication have been reported in \cite{Hafner2013} (using hybrid system theory) and in \cite{Englund2016a} where the authors pursue a virtual platooning approach at the Grand Cooperative Driving Challenge (GCDC) 2016. To the authors' best knowledge, only few experiments on \mbox{(semi-)}distributed MPC schemes for intersection automation, which rely on V2V communication, have been carried out, see e.g. \cite{Hult2019a}. 
In that regard, our contribution can be stated as: 
\begin{itemize}
	\item We propose and experimentally evaluate a fully distributed MPC scheme, which is independent of a central node as opposed to semi-distributed concepts \cite{Hult2019a}. 
\end{itemize}
That way, we contribute with a novel and relevant perspective to the sparse literature on experiments in that area of research.} To keep complexity at a manageable level, similar to \cite{Hult2019a}, we focus on scenarios in which agents cross the intersection straight. We then only manipulate the longitudinal acceleration while steering control can be taken care of by the driver. \chgRevIV{That said, the control system can be viewed as an adaptive cruise control (ACC) system which accommodates crossing vehicles.}

The remainder of the article is organized as follows. \prettyref{sec:problem} defines the intersection coordination problem along with a control-oriented kinematic agent model. 
Starting with a centralized problem formulation in \prettyref{sec:cmpc}, we continue with a distributed version and its fast numerical solution in \prettyref{sec:dmpc}. Thereafter, \prettyref{sec:archictecture} outlines the control system architecture which is utilized in our experiments and involves localization, communication and optimal control. Experimental results are finally discussed in \prettyref{sec:expResults}.

\section{Intersection Automation Problem} 
\label{sec:problem}

\subsection{\chgRevI{Notation}}
\label{sec:problem_notation}
\chgRevI{With $x_{k+j\mid k}$, we refer to the prediction of variable $x$ at the future time step 
$k+j$ given information up to time $k$. 
For $x\in\mathbb{R}^n$ and $i\in\{1,\ldots, n\}$, $[x]_i$ is the $i$-th entry of $x$, and the interval $[a,b] \subset \mathbb{N}$ with $a<b$ is denoted as $\mathbb{N}_{[a,b]}$.
Moreover, $\mathbb{N}^+$ is the set of positive integers and  \(A^{\T}\) denotes the transpose of a matrix \(A{}\in{}\mathbb{R}^{m\times n}\).}

\subsection{Problem Description}
\label{sec:problem_description}

\begin{problem}
We aim to automate agents in a four way, single lane unsignalized intersection by manipulating their acceleration. With fully automated longitudinal control, the agents shall cross the intersection \chgRevIV{straight} without any collisions while tracking a desired speed as close as possible. 
\end{problem}	
In the remainder, we rely on the following assumptions.
\begin{assumption}
	A1. Only single intersection scenarios \chgRevI{with one lane and one agent per direction are considered}; \chgAlex{A2.} The control scheme only influences longitudinal control to avoid collisions; \chgAlex{A3.} Lateral vehicle control is accommodated by a separate, independent control module (or the driver); \mbox{\chgAlex{A4.} The} desired route of every agent is determined by a high-level route planning algorithm; \chgAlex{A5.} Every agent is equipped with V2V communication; \chgAlex{A6.} Communication failures or package dropouts are neglected;  \chgAlex{A7.} The local MPC solutions at time \chgRevI{step} $k$ are available to all agents at time step $k+1$; {\chgAlex{A8.} The local MPCs are executed synchronously.} \chgAlex{A9.} Vehicle states are measurable or can be estimated appropriately. \chgAlex{A10. Avoiding rear-end collisions with frontal vehicles is not in scope.}
\end{assumption}

Assumptions A1-\chgAlex{A3, A5, A6 and A10} are common in the literature and are used to reduce complexity~\cite{Campos2014,Hult2016}. The use of a high-level planning algorithm which is postulated in \chgAlex{A4} is quite common in AV architectures too \cite{Lim2018}. Lastly, \chgAlex{A7} can be satisfied by choosing the MPC sampling time appropriately and \chgAlex{A8-A9} can be accomplished as shown in \prettyref{sec:archictecture}.

\vspace*{-2mm}
\subsection{Modeling of Agent Kinematics}
\label{sec:problem_agentKinematics}

\begin{figure}[t!]	
	\begin{center}
		\def\svgwidth{8.0cm}	
		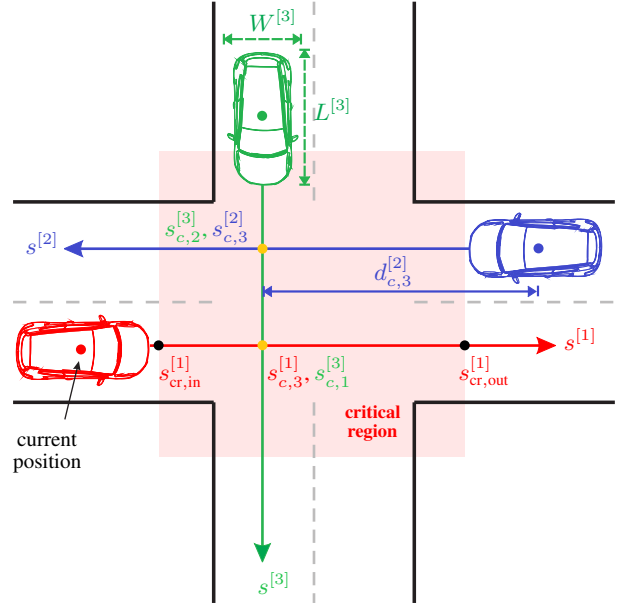
		\vspace*{-3mm}
		\caption{Example intersection scenario with $N_A=3$ \chgRevI{straight crossing} agents. The collision point with Agent $l$ along the path of Agent $i$ \chgAlex{within the critical region $[s_{\text{cr,in}}^{[i]},s_{\text{cr,out}}^{[i]}]$ is denoted as $s_{c,l}^{[i]}$ while $d_{c,l}^{[i]}$ is the distance to $s_{c,l}^{[i]}$}}\vspace*{-7mm} 
\label{fig:problem_agentKinematics_scheme}		
	\end{center}
\end{figure}

\chgAlex{For intersection automation, we consider the set \mbox{$\mathcal{A} \triangleq \{1, \ldots , N_A \}$} of $N_A$ connected and automated agents.}  
The motion dynamics of every agent $i$ is described \chgAlex{in terms of its geometric center's} acceleration $a_{x}^{[i]}$, velocity $v^{[i]}$ and path coordinate $s^{[i]}$ along a given path, see \prettyref{fig:problem_agentKinematics_scheme}. 
For these kind of problems, it is a common approach to describe the time evolution of velocity and position as a double integrator \cite{Campos2014,Hult2016,Molinari2018}. By modeling drivetrain dynamics as a first order lag element, Agent $i$'s motion can be summarized as a linear time-invariant state space model, i.e.,
\begin{align}
\frac{d}{dt}
\begin{bmatrix}
{a}_{x}^{[i]} \\
{v}^{[i]} \\
{s}^{[i]}
\end{bmatrix} &= 
\underbrace{\begin{bmatrix}
	-\frac{1}{T_{a_x}^{[i]}} & 0 & 0 \\
	1 & 0 & 0 \\
	0 & 1 & 0
	\end{bmatrix}}_{A^{[i]}}
\underbrace{
	\vphantom{\begin{bmatrix}
		\frac{1}{T_{a_x}^{[i]}} \\
		0 \\
		0
		\end{bmatrix}}
	\begin{bmatrix}
	{a}_{x}^{[i]} \\
	{v}^{[i]} \\
	{s}^{[i]} 
	\end{bmatrix}}_{x^{[i]}} +
\underbrace{\begin{bmatrix}
	\frac{1}{T_{a_x}^{[i]}} \\
	0 \\
	0
	\end{bmatrix}}_{B^{[i]}} \underbrace{\vphantom{\begin{bmatrix}
		\frac{1}{T_{a_x}^{[i]}} \\
		0 \\
		0
		\end{bmatrix}}
	a_{x,\text{ref}}^{[i]}}_{u^{[i]}}
\label{eq:problem_agentKinematics_ssContinuous}
\end{align}
where 
\(
x^{[i]} \triangleq
[{a}_{x}^{[i]},\,
{v}^{[i]},\,
{s}^{[i]}]^{\T}
\)
is the state vector, 
\mbox{$u^{[i]}\triangleq a_{x,\text{ref}}^{[i]}$}
the reference acceleration as control input and $T_{a_{x}}^{[i]}$ the dynamic drivetrain time constant. States and inputs are constrained by polyhedral sets,
that is, 
\(
x^{[i]} \in \mathcal{X}^{[i]} \subseteq \mathbb{R}^{n_x}
\)
and
\mbox{$u^{[i]} \in \mathcal{U}^{[i]} \subseteq \mathbb{R}^{n_u}$}
with $n_x = 3$ and $n_u = 1$. To be used within numerical optimization algorithms, we discretize \eqref{eq:problem_agentKinematics_ssContinuous} by using zero-order hold discretization. This way, we gain the discrete-time linear time-invariant state space model 
\begin{align}
x_{k+1}^{[i]} = A_d^{[i]} x_{k}^{[i]} + B_d^{[i]} u_{k}^{[i]}
\label{eq:problem_agentKinematics_ssDiscrete}
\end{align}
with
\(
A_d^{[i]} \triangleq e^{A^{[i]}T_s}
\)
and 
\(
B_d^{[i]} \triangleq \int_0^{T_s} e^{A^{[i]} \tau} d\tau B^{[i]} 
\) where $T_s > 0$ is the corresponding sample time.

\subsection{Distance Between Agents}
\label{sec:problem_agentDistance}

To compute the distance between two agents $i,l \in \mathcal{A}$, in a first step, the collision points $s_{c,l}^{[i]}$ and $s_{c,i}^{[l]}$ have to be determined. 
According to \prettyref{fig:problem_agentKinematics_scheme}, these collision points correspond to the intersection of the agents' paths along their path coordinates $s^{[i]}$ and $s^{[l]}$. If their respective paths do not intersect, we define $s_{c,l}^{[i]} = s_{c,i}^{[l]} = \infty$. In a second step, agents $i$ and $l$ calculate the distances $d_{c,l}^{[i]}$ and $d_{c,i}^{[l]}$ to their respective collision points $s_{c,l}^{[i]}$ and $s_{c,i}^{[l]}$ respectively, that is,
\begin{align}
d_{c,l}^{[i]} \triangleq \begin{cases} \lvert s^{[i]} - s_{c,l}^{[i]} \rvert & ,s_{c,l}^{[i]} \neq \infty \\
\infty & ,\text{otherwise}. \end{cases} \label{eq:problem_agentDistance_distAgentToCP}
\end{align}
Finally, we define the distance between Agent $i$ and Agent $l$ as the sum of distances to their joint collision point, i.e., 
\begin{align}
\mathrm{dist}(i,l) \triangleq d_{c,l}^{[i]} + d_{c,i}^{[l]}.
\label{eq:problem_agentDistance_distAgentToAgent}
\end{align}

\section{Centralized Problem Formulation}
\label{sec:cmpc}

In this section, we aim to formalize the problem definition in \prettyref{sec:problem_description} in terms of a centralized OCP, that is, a problem that is solved in a receding horizon fashion by a central node. For reasons outlined in \prettyref{sec:introduction_contribution}, we rely on an MPC-based framework. When applying MPC, at every time step $k$ we solve a finite-time OCP over a prediction horizon of $N$ time steps. After optimization, only the first control input is applied to the plant and optimization is repeatedly executed over a shifted horizon at time $k+1$. With a centralized OCP in place, we transition to a fully distributed formulation in \prettyref{sec:dmpc}.

\subsection{Agent Objectives and Constraints}
\label{sec:cmpc_objCons}

The centralized OCP is actually an aggregation of every agent's local objectives and constraints as well as joint collision avoidance constraints which couple the agents. 

\vspace*{0.5mm}
{\textbf{Objectives}~~In terms of local objectives, every agent $i$ is intended to follow a reference speed $v_{\text{ref}}^{[i]}$ while at the same time fertilizing ride comfort and efficiency by minimizing step changes of the control input (i.e., the longitudinal acceleration $u^{[i]}=a_x^{[i]}$) and its magnitude, respectively. These objectives can be cast as a convex quadratic cost of the form} 
\begin{align} \label{eq:cmpc_objCons_costFcn}
J^{[i]}(x_{\cdot \mid k}^{[i]}, u_{\cdot \mid k}^{[i]}) &
\triangleq \,Q_N^{[i]} \,(v_{\text{ref},k+N\mid k}^{[i]} - v_{k+N\mid k}^{[i]})^2 \notag\\
&+ \,Q^{[i]} \sum^{N-1}_{j=1} (v_{\text{ref},k+j\mid k}^{[i]} - v_{k+j\mid k}^{[i]})^2 \\
&+ \,R^{[i]} \sum^{N-1}_{j=0} ( \Delta u_{k+j\mid k}^{[i]})^2 + S^{[i]} \sum^{N-1}_{j=0} (u_{k+j\mid k}^{[i]})^2 \notag
\end{align}
where the first term represents the terminal cost, $Q^{[i]}$, $Q_N^{[i]}$, $R^{[i]}$, $S^{[i]} > 0$ are positive scalar weights \chgRevI{and 
\mbox{\(
\Delta u_{k+j\mid k}^{[i]} \triangleq u_{k+j\mid k}^{[i]} - u_{k+j-1\mid k}^{[i]}
\)}
is the step change of control inputs with $u_{k-1\mid k}^{[i]} \triangleq u_{k-1}^{[i]}$ for $j=0$.}

\vspace*{0.5mm}
{\textbf{Constraints}~~Besides objectives, we also need to accommodate constraints on the agents' inputs and states. Particularly, we contemplate actuator limitations in terms of maximum and minimum feasible accelerations which translates into box constraints on the inputs, that is,
\begin{align} 
\underline{u}^{[i]} \leq u_{k+j\mid k} \leq \overline{u}^{[i]},~~   \forall j \in \mathbb{N}_{[0,N-1]}. \label{eq:cmpc_objCons_inputCons}
\end{align} 
Moreover, it is intended to solely drive in the forward direction and to accommodate a maximum speed $\overline{v}^{[i]}$ (e.g., the road speed limit).
We phrase these conditions as a state constraint on the velocity, i.e.,
\begin{align} 
0 \leq v_{k+j\mid k}^{[i]} \leq \overline{v}_{k+j\mid k}^{[i]},~~ \forall j \in \mathbb{N}_{[1,N]}. \label{eq:cmpc_objCons_stateCons}
\end{align}
\chgRevI{To guarantee collision avoidance, we need to satisfy the following condition before Agent $i$ enters the critical region $[s_{\text{cr,in}}^{[i]},s_{\text{cr,out}}^{[i]}]$ of the intersection (see \prettyref{fig:problem_agentKinematics_scheme}): At the end of the prediction horizon, Agent $i$ has either i) left the critical region or ii) has stopped before entering that region \cite{Katriniok2017a}. \mbox{Case i)}} is equivalent to a terminal constraint which forces Agent $i$ to leave the critical region at time step $k+N$, that is, 
\begin{align} 
s_{k+N\mid k}^{[i]} \geq s_{\text{cr,out}}^{[i]}.
\label{eq:cmpc_objCons_vMinMeanCons}
\end{align}
That said, we only have to impose \eqref{eq:cmpc_objCons_vMinMeanCons} if Agent $i$ is about to enter or is located within the critical region, 
that is, if \mbox{$s_k^{[i]} \in [ s_{c,\text{in}}^{[i]} - d_{\text{brake}}, \,s_{c,\text{out}}^{[i]}]$} where $d_{\text{brake}}>0$ is a brake safe distance. In all other cases, i.e., far away from the intersection and after crossing it, we set $s_{\text{cr,out}}^{[i]}$ in \eqref{eq:cmpc_objCons_vMinMeanCons} to a sufficiently large negative value, thus satisfying the constraint at all times.
If constraint \eqref{eq:cmpc_objCons_vMinMeanCons} causes the OCP to be infeasible, the respective agent is forced to stop before entering the critical region \chgRevI{--- thus accommodating case ii)} and recovering feasibility. Finally, constraints \eqref{eq:cmpc_objCons_inputCons}, \eqref{eq:cmpc_objCons_stateCons} and \eqref{eq:cmpc_objCons_vMinMeanCons} can concisely be \mbox{written as}
\begin{align}
P_x^{[i]} x_{\cdot\mid k}^{[i]} + P_u^{[i]} u_{\cdot\mid k}^{[i]} + q_{xu}^{[i]}\leq 0
\label{eq:cmpc_objCons_stateInputConsCombined}
\end{align}
where $P_x^{[i]} \succeq 0$, $P_{u}^{[i]} \succeq 0$, and $q_{xu}^{[i]}$ are matrices and vectors of appropriate dimension.

\subsection{Collision Avoidance}
\label{sec:cmpc_CA}

While all constraints in \prettyref{sec:cmpc_objCons} refer to the individual Agent $i$, collision avoidance eventually couples the agents among each other. To mathematically claim collision avoidance, we first define Agent $i$'s conflict set 
\begin{align*}
\mathcal{A}_{c}^{[i]} \triangleq \Bigl\{  l \in \mathcal{A} \mid l \neq i \land s_{c,l}^{[i]} \neq \infty \Bigr\}, 
\end{align*}
i.e., the set of agents $l\neq i$ which have a joint collision point with Agent $i$.
Avoiding collisions between Agent $i$ and \mbox{Agent $l \in \mathcal{A}_c^{[i]}$} is then stated as a lower bound on their distance
\begin{align}
d_{c,l,k+j\mid k}^{[i]} + d_{c,i,k+j\mid k}^{[l]} \geq d_{\text{safe}},~ \forall j \in \mathbb{N}_{[1,N]} 
\label{eq:cmpc_CA_CACons}
\end{align}
where $d_{c,l,k+j\mid k}^{[i]}$, $d_{c,i,k+j\mid k}^{[l]}$ at the predicted time step $k+j$ depend on the predicted path coordinates $s_{k+j\mid k}^{[i]}$, $s_{k+j\mid k}^{[l]}$ and the collision points $s_{c,l}^{[i]}$, $s_{c,i}^{[l]}$ 
in accordance to \eqref{eq:problem_agentDistance_distAgentToCP}. 
Moreover, $d_{\text{safe}} > 0$ is a suitable safety distance which also captures the agents' width and length.

\subsection{Centralized Optimal Control Problem}
\label{sec:cmpc_centralOCP}

The centralized intersection coordination problem results from the sum of agents' costs \eqref{eq:cmpc_objCons_costFcn} subject to their state and input constraints \eqref{eq:cmpc_objCons_stateInputConsCombined}, their dynamics \eqref{eq:problem_agentKinematics_ssDiscrete} and their joint collision avoidance constraints \eqref{eq:cmpc_CA_CACons}. In essence, we obtain
\begin{subequations} \label{eq:cmpc_centralOCP_defOCP}
	\begin{align}
	\hspace*{-3.4mm}\underset{ x_{\cdot\mid k},\, u_{\cdot\mid k} }{\mathmin} ~&  \sum_{i=1}^{N_A} J^{[i]}({x}_{\cdot \mid k}^{[i]}, {u}_{\cdot \mid k}^{[i]};\, x_k^{[i]}) \label{eq:cmpc_centralOCP_defOCP_cost}\\ 
	\hspace*{-3.4mm}\mathst~& \textbf{agent constraints} \text{ -- }  \forall i \in \mathcal{A}: \notag\\
	\hspace*{-3.4mm}& P_x^{[i]} x_{\cdot\mid k}^{[i]} + P_u^{[i]} u_{\cdot\mid k}^{[i]} + q_{xu}^{[i]}\leq 0 \label{eq:cmpc_centralOCP_defOCP_stateInputCons} \\
	\hspace*{-3.4mm}~& x_{k+j+1\mid k}^{[i]} \hspace*{-0.8mm}=\hspace*{-0.8mm} A_d^{[i]} x_{k+j\mid k}^{[i]} \hspace*{-0.8mm}+\hspace*{-0.8mm} B_d u_{k+j\mid k}^{[i]},\forall j \hspace*{-0.5mm} \in \hspace*{-0.5mm} \mathbb{N}_{[0,N-1]}
	\label{eq:cmpc_centralOCP_defOCP_sysDynCons}\\
	\hspace*{-3.4mm}~& x_{k\mid k}^{[i]} = x_k^{[i]} \label{eq:eq:cmpc_centralOCP_defOCP_x0}\\
	\hspace*{-3.4mm} ~& \textbf{coupling constraints} \text{ -- }  \forall i \in \mathcal{A},~ \forall l \in {\mathcal{A}}_{c}^{[i]}:\notag\\
	\hspace*{-3.4mm}& d_{c,l,k+j\mid k}^{[i]} + d_{c,i,k+j\mid k}^{[l]} \geq d_{\text{safe}},~\forall j \in \mathbb{N}_{[1,N]}
	\label{eq:cmpc_centralOCP_defOCP_CAcons}
	\end{align}
\end{subequations}
In \eqref{eq:cmpc_centralOCP_defOCP}, the quadratic cost \eqref{eq:cmpc_centralOCP_defOCP_cost} is convex as all weights are positive. The same holds for the agent constraints as these are linear in the decision variables. Eventually, nonconvexity arises with the absolute value collision avoidance constraint \eqref{eq:cmpc_centralOCP_defOCP_CAcons} which relates 
to every agent's decision to cross the intersection before or after the other agent. To ensure feasibility of \eqref{eq:cmpc_centralOCP_defOCP}, given that terminal constraint \eqref{eq:cmpc_objCons_vMinMeanCons} is satisfied, constraints \eqref{eq:cmpc_centralOCP_defOCP_stateInputCons} and \eqref{eq:cmpc_centralOCP_defOCP_CAcons} are implemented as soft constraints.

\section{Fully Distributed MPC Scheme}
\label{sec:dmpc}

\subsection{Decomposition}
\label{sec:dmpc_decomposition}

\subsubsection{Decoupling Collision Avoidance Constraints}
\label{sec:dmpc_decomposition_CAdecouple}
When distributing problem \eqref{eq:cmpc_centralOCP_defOCP}, it would be detrimental to minimize every agent's cost subject to its input and state constraints, and to impose  collision avoidance constraints \eqref{eq:cmpc_centralOCP_defOCP_CAcons} on \mbox{Agent $i$} and Agent $l$ simultaneously when their paths intersect. In that case, both agents have to take independent decisions at time $k$ instead of  jointly \textit{agreeing} on how to avoid collisions --- this may eventually lead to collisions.

To accommodate this issue, we define the bijective prioritization function $\gamma: \mathcal{A} \rightarrow \mathcal{A}$ which assigns a unique priority to every agent where a lower value corresponds to a higher priority \cite{Katriniok2017a}. To this end, we specify the prioritized conflict set 
\begin{align*}
\mathcal{A}_{c,\gamma}^{[i]} \triangleq \bigl\{  l \in {\mathcal{A}}_{c}^{[i]} \mid \gamma(l) < \gamma(i) \bigr\}
\end{align*}
containing all agents $l \in \mathcal{A}_c^{[i]}$ which have a joint collision point with Agent $i$ but a higher priority. For all agents $l \in \mathcal{A}_{c,\gamma}^{[i]}$, we then impose constraint \eqref{eq:cmpc_centralOCP_defOCP_CAcons} only on Agent $i$ (having lower priority), thus yielding fully decoupled agent OCPs.

\begin{remark}
We would like to stress that agent priorities do not imply an intersection crossing order. Conversely, they just define which agent has to accommodate collision avoidance constraints. 
\end{remark}

\subsubsection{Reformulation of Collision Avoidance Constraints}
\label{sec:dmpc_decomposition_CAreform}

To appropriately pose collision avoidance constraint \eqref{eq:cmpc_centralOCP_defOCP_CAcons} for numerical optimization, we rewrite it in a quadratic form. After rearranging  \eqref{eq:cmpc_centralOCP_defOCP_CAcons}, we square both sides of the inequality, i.e.,
\begin{align}
(s_{k+j\mid k}^{[i]} - s_{c,l}^{[i]})^2 \geq (d_{\text{safe}} - d_{c,i,k+j\mid k}^{[l]})^2. \label{eq:dmpc_decomposition_CAreform_sqCAcons}
\end{align}
It can be recognized that applying \eqref{eq:dmpc_decomposition_CAreform_sqCAcons} as constraint is only valid if $d_{\text{safe}} - d_{c,i,k+j\mid k}^{[l]} \geq 0$ holds. Conversely, \eqref{eq:dmpc_decomposition_CAreform_sqCAcons} only needs to be imposed if $d_{\text{safe}} - d_{c,i,k+j\mid k}^{[l]} > 0$. In all other cases, the original constraint \eqref{eq:cmpc_centralOCP_defOCP_CAcons} is satisfied \textit{per se}. That said, we define the time dependent prioritized conflict set
\begin{align*}
{\mathcal{A}}_{c,\gamma,j}^{[i]} \triangleq \Bigl\{  l \in \mathcal{A}_{c,\gamma}^{[i]} \mid d_{\text{safe}} - d_{c,i,k+j\mid k}^{[l]} > 0 \Bigr\} 
\end{align*}
and impose \eqref{eq:dmpc_decomposition_CAreform_sqCAcons} for  every $l \in {\mathcal{A}}_{c,\gamma,j}^{[i]}$ and \mbox{$j \in \mathbb{N}_{[1,N]}$} to ensure collision avoidance. Moreover, we specify the parameter sequence 
\mbox{
	\(
	z_{\cdot \mid k}^{[i]}\triangleq \{ (d_{c,i,k+j\mid k}^{[l]} )_{j\in\mathbb{N}_{[1,N]}} \}_{{l \in {\mathcal{A}}_{c}^{[i]}}}
	\)}, 
which \mbox{Agent $i$} receives from the other agents via V2V communication, and rewrite \eqref{eq:dmpc_decomposition_CAreform_sqCAcons} in dependence of Agent $i$'s state vector $x_{k+j\mid k}^{[i]}$ as
\begin{align}
(x_{k+j\mid k}^{[i]})^{\T} P_{l,j}^{[i]} \,x_{k+j\mid k}^{[i]} + (q_{l,j}^{[i]})^{\T} 
\,x_{k+j\mid k}^{[i]} + r_{l,j}^{[i]}(z_{\cdot\mid k}^{[i]}) \leq 0 
\label{eq:dmpc_decomposition_CAreform_sqCAconsMatVec}
\end{align}
with $P_{l,j}^{[i]} \triangleq \mathrm{diag}(0,0,-1)$, a suitable vector $q_{l,j}^{[i]} \in \mathbb{R}^{n_x}$ and parameterized scalar $r_{l,j}^{[i]}(z_{\cdot\mid k}^{[i]}) \in \mathbb{R}$. 
Evidently, $P_{l,j}^{[i]}$ is a negative semi-definite matrix which reflects the  nonconvexity of the collision avoidance constraint.  
We iteratively substitute system dynamics
\begin{align}
\hspace*{-1mm}
x_{k+j\mid k}^{[i]}( {u}_{\cdot\mid k}^{[i]}  )
{}={} 
(A_d^{[i]})^{j} x_{k\mid k}^{[i]}
{}+{} 
\sum_{\iota{}={}0}^{j-1} (A_d^{[i]})^{j-1-\iota}B_{d}^{[i]} u_{k+\iota\mid k}^{[i]} 
\label{eq:dmpc_decomposition_CAreform_sysDynamics}
\end{align}
into \eqref{eq:dmpc_decomposition_CAreform_sqCAconsMatVec}, thus yielding
\begin{align}
(u_{\cdot\mid k}^{[i]})^{\T} \bar{P}_{l,j}^{[i]}\, u_{\cdot\mid k}^{[i]} + (\bar{q}_{l,j}^{[i]})^{\T} u_{\cdot\mid k}^{[i]} + \bar{r}_{l,j}^{[i]}(z_{\cdot\mid k}^{[i]}) \leq 0
\label{eq:dmpc_decomposition_CAreform_condensedCAcons}
\end{align}
with suitable matrices, vectors and scalars $\bar{P}_{l,j}^{[i]}\preceq0$, $\bar{q}_{l,j}^{[i]}$ and $\bar{r}_{l,j}^{[i]}$, where  $\bar{r}_{l,j}^{[i]}$ is parameterized with respect to $z_{\cdot\mid k}^{[i]}$.
\vspace*{1mm}

\subsubsection{Distributed Optimal Control Problem}
\label{sec:dmpc_decomposition_docp}

Finally, for every \mbox{Agent $i$} we substitute \eqref{eq:dmpc_decomposition_CAreform_sysDynamics} into the cost \eqref{eq:cmpc_centralOCP_defOCP_cost} as well as input and state constraints \eqref{eq:cmpc_centralOCP_defOCP_stateInputCons}, 
thus obtaining the condensed cost $\bar{J}^{[i]}$ in  \eqref{eq:dmpc_decomposition_docp_defDOCP_cost} and the condensed constraint \eqref{eq:dmpc_decomposition_docp_defDOCP_stateInputCons}. 
This way, the resulting OCP, which is solved in parallel by every \mbox{agent $i$}, can be stated as a parameterized nonconvex QCQP, i.e., 
\vspace*{-1mm}
\begin{subequations} \label{eq:dmpc_decomposition_docp_defDOCP}
	\begin{align}
	\hspace*{-2mm}\underset{ u^{[i]}_{\cdot\mid k}}{\mathmin} ~& ~  \bar{J}^{[i]}({u}_{\cdot \mid k}^{[i]};\, x_k^{[i]} ) \label{eq:dmpc_decomposition_docp_defDOCP_cost}\\ 	
	\hspace*{-2mm}\mathst~&~~~~~~~~~~~~~~~~~~~~~~~~\, \bar{P}_{xu}^{[i]}\, u_{\cdot\mid k}^{[i]} + \bar{q}_{xu}^{[i]} \leq 0 \label{eq:dmpc_decomposition_docp_defDOCP_stateInputCons}\\ 
	~&~(u_{\cdot\mid k}^{[i]})^{\T} \bar{P}_{l,j}^{[i]}\, u_{\cdot\mid k}^{[i]} + (\bar{q}_{l,j}^{[i]})^{\T} u_{\cdot\mid k}^{[i]} + \bar{r}_{l,j}^{[i]} \leq 0, \label{eq:dmpc_decomposition_docp_defDOCP_CAcons}\\[-1mm]
	~&~ ~~~~~~~~~~~~~~~~~\, \forall l \in \mathcal{A}_{c,\gamma,j}^{[i]},~  \forall j \in \mathbb{N}_{[1,N]} \notag
	\end{align}
\end{subequations}
where $\bar{P}_{xu}^{[i]} \succeq 0$ and $\bar{q}_{xu}^{[i]}$ are matrices and vectors of appropriate dimension. To ensure the feasibility of OCP \eqref{eq:dmpc_decomposition_docp_defDOCP}, given that terminal constraint \eqref{eq:cmpc_objCons_vMinMeanCons} is satisfied (see \prettyref{sec:cmpc_objCons}), state constraints \eqref{eq:dmpc_decomposition_docp_defDOCP_stateInputCons} and collision avoidance constraints \eqref{eq:dmpc_decomposition_docp_defDOCP_CAcons} are implemented as soft constraints. 

\begin{remark}
	Along the lines of agent prioritization, in the distributed setting terminal constraint \eqref{eq:cmpc_objCons_vMinMeanCons} has only to be imposed if there are conflicting agents of higher priority that have not yet left the critical region at time $k$.
\end{remark}

\subsection{Fast Numerical Solution of the Distributed OCP}
\label{sec:dmpc_solution}

Referring to \prettyref{sec:introduction_contribution}, a major challenge is to solve the nonconvex problem \eqref{eq:dmpc_decomposition_docp_defDOCP} fast on an embedded hardware. To this end, we rely on the penalty convex-concave procedure (CCP) \cite{Lipp2016} which allows us to solve \eqref{eq:dmpc_decomposition_docp_defDOCP} as a sequence of QPs and as such to leverage 
efficient of-the-shelve QP solvers.

\subsubsection{Background on the Convex-Concave Procedure} 
\label{sec:dmpc_solution_CCPbackgrnd}
The main idea behind CCP is that every nonconvex function \mbox{$h: \mathbb{R}^n \rightarrow \mathbb{R}$} can be written as the difference of two convex functions $f: \mathbb{R}^n \rightarrow \mathbb{R}$ and \mbox{$g: \mathbb{R}^n \rightarrow \mathbb{R}$}, that is,
\vspace*{-0.5mm}
\begin{align*}
	h(x) = f(x) - g(x).
\end{align*}
Along these lines, a nonconvex minimization problem with a smooth, nonconvex cost $h_0(x)$ and smooth, nonconvex constraints $h_i(x) \leq 0$ with $i \in \mathbb{N}_{[1,M]}$ can be written as 
\begin{subequations} \label{eq:dmpc_solution_CCPbackgrnd_optCCPidea}
	\begin{align}
	\hspace*{-2mm}\underset{ x }{\mathmin} ~& ~  f_0(x) - g_0(x) \\ 
	\hspace*{-2mm}\mathst~&~ f_i(x) - g_i(x) \leq 0, ~~ \forall i \in \mathbb{N}_{[1,M]}.
	\end{align}
\end{subequations}
Then, starting at an initial point $x^0$, optimization problem \eqref{eq:dmpc_solution_CCPbackgrnd_optCCPidea} is solved iteratively by successively linearizing the nonconvex term $-g_i(x)$, $i \in \mathbb{N}_{[0,M]}$ with respect to the current solution candidate $x^\nu$ at every iteration $\nu$, that is,
\begin{align}
\tilde{h}_i(x) \triangleq f_i(x) - (\nabla_{\hspace*{-0.5mm}x} \,g_i(x^\nu))^{\T}  x - g_i(x^\nu), ~ i \in \mathbb{N}_{[0,M]}
\label{eq:dmpc_solution_convexifiedFcn}
\end{align}
until a suitable convergence criterion is satisfied. The penalty CCP method, applied in the remainder, is a variant of the standard CCP method and does not require a feasible initial point. As such a point may not even be known \textit{a priori}, penalty CCP is much better suited for our type of application.

\subsubsection{Application of the Penalty CCP Method} 
\label{sec:dmpc_solution_CCPappl}

To solve the nonconvex \mbox{OCP \eqref{eq:dmpc_decomposition_docp_defDOCP}} fast, we apply the penalty CCP method to approach a numerical solution in an iterative fashion, see \prettyref{algo:dmpc_solution_CCPappl_CCPalgo}.
For reasons of clarity, the superscript $[i]$, indicating the corresponding agent, is omitted in \prettyref{algo:dmpc_solution_CCPappl_CCPalgo} and the subsequent description of the algorithm. 
At every time step $k$, every agent $i \in \mathcal{A}$ starts with an initial solution candidate $u^0 \triangleq u_{\cdot\mid k}^0$ (representing the control actions over the prediction horizon), an initial penalty weight $\rho_c^0>0$, a maximum penalty $\overline{\rho}_c\geq \rho_c^0$ and a penalty update parameter $\mu>1$. These parameters are further explained in the remainder of the section. That said, the following steps are carried out.

\textbf{Step 1~~} 
In problem \eqref{eq:dmpc_decomposition_docp_defDOCP}, the cost \eqref{eq:dmpc_decomposition_docp_defDOCP_cost} as well as input and state constraints \eqref{eq:dmpc_decomposition_docp_defDOCP_stateInputCons} are convex while only collision avoidance constraints \eqref{eq:dmpc_decomposition_docp_defDOCP_CAcons} are nonconvex. By virtue of \prettyref{sec:dmpc_decomposition_CAreform}, the nonconvex part of \eqref{eq:dmpc_decomposition_docp_defDOCP_CAcons} corresponds to the term $(u_{\cdot\mid k})^{\T} \bar{P}_{l,j} \, u_{\cdot\mid k}$ with $\bar{P}_{l,j} \preceq 0$. By means of \eqref{eq:dmpc_solution_convexifiedFcn}, we have to linearize this term with respect to the current solution candidate $u^\nu$, thus gaining 
\(
(u^\nu)^{\T} \bar{P}_{l,j} \, u + (u^\nu)^{\T} \bar{P}_{l,j} \,u^\nu
\). 

\setlength{\textfloatsep}{4.5mm}
\begin{algorithm}[b!]
	\caption{Penalty CCP algorithm to solve every agent's OCP \eqref{eq:dmpc_decomposition_docp_defDOCP}. For reason of clarity, we drop the superscript $[i]$.} 
	\begin{algorithmic}
		\State For every agent, we run the following steps at time $k$:
		\State \textbf{Input:} Initial point $u^0 \triangleq u^0_{\cdot\mid k}$, $\rho_c^0 > 0$, $\overline{\rho}_c$, and $\mu>1$ \vspace*{1.5mm}
		\State Set iteration $\nu \leftarrow 0$ and define $u^\nu \triangleq u^\nu_{\cdot\mid k}$.
		\Repeat	\vspace*{0.5mm}
		\State \textbf{Step 1.} Convexify collision avoidance constraints \eqref{eq:dmpc_decomposition_docp_defDOCP_CAcons}:
		\vspace*{-2mm} 
		\begin{equation}
		\left( (u^\nu)^{\T} \bar{P}_{l,j} + \bar{q}_{l,j}^{\T}\right) u + \bar{r}_{l,j} + {(u^\nu)^{\T} \bar{P}_{l,j} \,u^\nu} \leq 0
		\label{eq:dmpc_solution_CCPappl_CCPalgo_convexCAcons}
		\end{equation}
		\State \textbf{Step 2.} Solve convexified version of problem \eqref{eq:dmpc_decomposition_docp_defDOCP}: \vspace*{1.5mm}
		\State $(u^\nu, \epsilon_{x}^\nu, \epsilon_{c}^\nu) \leftarrow$
		\vspace*{-4mm}
		\begin{subequations} \label{eq:dmpc_solution_CCPappl_CCPalgo_convexOCP}
		\begin{align}
		~~~\text{arg\,} \hspace*{-1.5mm}\underset{ u, \epsilon_{x}, \epsilon_{c} }{\mathmin} &   {\bar{J}}({u};  x_k ) + \rho_{x} \epsilon_{x} + \rho_c^\nu \sum_{j=1}^{N} [\epsilon_{c}]_j \label{eq:dmpc_solution_CCPappl_CCPalgo_convexOCP_cost} \\
		~~~\mathst~	& \chgRevIV{\textbf{agent state \& input constraints:}}\notag\\
		& \chgRevIV{\bar{P}_{xu} u + \bar{q}_{xu} \leq {\bar\Xi_{x}\epsilon_{x}},~~\epsilon_{x} \geq 0} \\[1mm]
		& \chgRevIV{\textbf{collision avoidance \eqref{eq:dmpc_solution_CCPappl_CCPalgo_convexCAcons}: }\hspace*{-0.5mm} \forall l \hspace*{-0.5mm} \in \hspace*{-0.5mm} \mathcal{A}_{c,\gamma,j},\hspace*{0.1mm} \forall j \hspace*{-0.5mm} \in \hspace*{-0.5mm} \mathbb{N}_{[1,N]}} \notag\\ 
		& \hspace*{-4.3mm}\chgRevIV{\left( (u^\nu)^{\T} \bar{P}_{l,j} + \bar{q}_{l,j}^{\T}\right) u + \bar{r}_{l,j} + (u^\nu)^{\T} \bar{P}_{l,j} \,u^\nu \leq [\epsilon_c]_j,} \notag\\
		& \chgRevIV{\epsilon_{c} \geq 0}
		\end{align}
		\end{subequations}
		\State \textbf{Step 3.} Update weight: ~~~~~~~~~~$\rho_c^{\nu+1} \leftarrow \min\{\mu\rho_c^\nu,\overline{\rho}_c\}$
		\State ~~~~~~~~~\,Update iteration count: ~~~~\,~$\nu \leftarrow \nu+1$ 
		\Until{stopping criterion is satisfied}\vspace*{1.5mm}
		\State \textbf{Output:} Stationary point $u^\star \leftarrow u^\nu$ 
	\end{algorithmic} 
	\label{algo:dmpc_solution_CCPappl_CCPalgo}
\end{algorithm}

\textbf{Step 2 \& Step 3~} 
By linearizing the nonconvex part of \eqref{eq:dmpc_decomposition_docp_defDOCP_CAcons}, the originally nonconvex constraints 
are replaced by linear (and as such convex) constraints \eqref{eq:dmpc_solution_CCPappl_CCPalgo_convexCAcons}. This way, the resulting OCP \eqref{eq:dmpc_solution_CCPappl_CCPalgo_convexOCP} is convex.  
In \eqref{eq:dmpc_solution_CCPappl_CCPalgo_convexOCP}, constraints \eqref{eq:dmpc_solution_CCPappl_CCPalgo_convexCAcons} are imposed as a soft constraints by introducing a vector of slack variables $\epsilon_c\geq0$ where each slack variable is penalized in the augmented cost function \eqref{eq:dmpc_solution_CCPappl_CCPalgo_convexOCP_cost}. 
In standard penalty CCP, every nonconvex constraint comes along with its own slack variable. 
To reduce the dimensionality and as such the computational complexity of OCP \eqref{eq:dmpc_solution_CCPappl_CCPalgo_convexOCP}, we propose to use a single slack variable per time step $k+j$.
That said, the corresponding slack variable represents the violation of \eqref{eq:dmpc_solution_CCPappl_CCPalgo_convexCAcons} with respect to the most critical agent (the $\infty$-norm) at that time step. 
Starting with a low initial penalty $\rho_c^0>0$, the numerical algorithm is allowed to initially explore potential local optima, while in subsequent iterations $\nu>0$ of the algorithm the penalty $\rho_c$ is increased up to a specified maximum $\overline{\rho}_c \geq \rho_c^0$ \mbox{({Step 3})}. That way, subsequent iterations \textit{tie} the solution to a certain (local) feasible region.  
By virtue of \cite{Lipp2016}, it can be shown that as $\epsilon_c \rightarrow 0$ the solution of OCP \eqref{eq:dmpc_solution_CCPappl_CCPalgo_convexOCP} provides a feasible point for the original nonconvex QCQP \eqref{eq:dmpc_decomposition_docp_defDOCP}. It is evident, though, that \prettyref{algo:dmpc_solution_CCPappl_CCPalgo} is a \textit{local method} to solve problem \eqref{eq:dmpc_decomposition_docp_defDOCP}. However, solving for a local stationary point instead of the global optimal solution is a common approach, especially in a \chgRevIV{real-time} setting. Besides $\epsilon_{c}$, we apply the additional slack variable $\epsilon_{x} \geq 0$ ($\infty$-norm, weighted by $\rho_{x}>0$ in the cost) along with a suitable vector $\bar{\Xi}_{x}$ to implement state \mbox{constraint \eqref{eq:cmpc_objCons_stateCons}} on the agent's velocity as a \mbox{soft constraint.}

\textbf{Stopping Criterion~~} The algorithm is iterated until a stopping criterion is satisfied. Particularly, we rely on the condition proposed in \cite{Lipp2016}, which holds if either the penalty $\rho$ has reached its maximum value $\overline{\rho}$ or the improvement of the objective function \eqref{eq:dmpc_solution_CCPappl_CCPalgo_convexOCP_cost} between iteration $\nu-1$ and $\nu$ is less than a sufficiently small threshold and 
the overall violation $\rho\, \|\epsilon_c^{\nu}\|_1$ of collision avoidance constraints is sufficiently small. 

After \prettyref{algo:dmpc_solution_CCPappl_CCPalgo} has converged to a stationary point $u_{\cdot\mid k}^\star$ at time $k$, only the first control action $u_{k\mid k}^\star$ is applied to the system. At the next time step $k+1$, \prettyref{algo:dmpc_solution_CCPappl_CCPalgo} is warm started by leveraging the shifted solution $u_{\cdot\mid k}^\star$ from time \mbox{step $k$.} 
A proof showing that \prettyref{algo:dmpc_solution_CCPappl_CCPalgo} finally converges to a (local) stationary point can be found in \cite{Lipp2016}.

\begin{remark} The advantage of (penalty) CCP over other methods, such as SQP, is that more information is retained in each of its iterates \cite{Lipp2016}. While SQP methods linearize the entire constraint \eqref{eq:dmpc_decomposition_docp_defDOCP_CAcons} at the current solution candidate, CCP only linearizes the nonconvex part $u^{\T} \bar{P}_{l,j} u$, thus preserving $\bar{q}_{l,j}^{\T} u + \bar{r}_{l,j}$ in \eqref{eq:dmpc_decomposition_docp_defDOCP_CAcons}.
\end{remark}

\section{In-Vehicle Control System Architecture}
\label{sec:archictecture}

\begin{figure*}[h!]
	\centering
	\def\svgwidth{180mm}
	\vspace*{-10mm}
	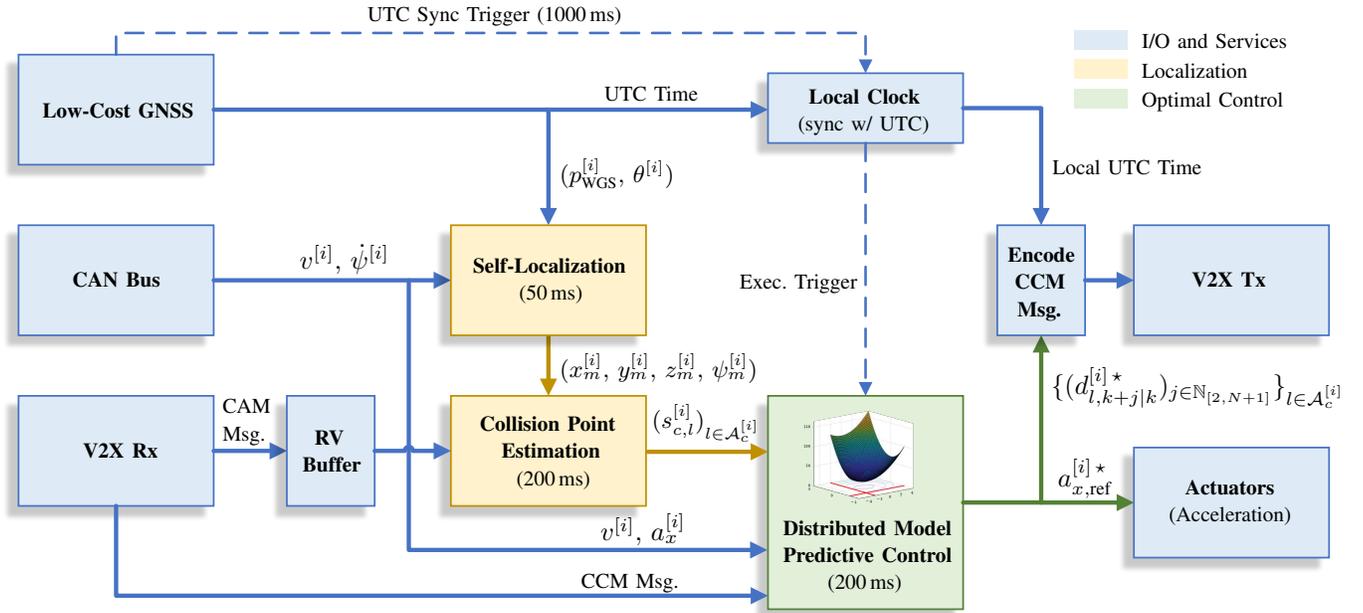	
	\vspace*{-3mm}
	\caption{Control system architecture utilized in every vehicle and run on a dSPACE MicroAutoBox II (illustration for Agent $i$). First, the input data is processed within localization algorithms before it is leveraged by the local MPC controller. The execution of the MPC is synchronized with the other agents. Finally, the optimized acceleration is applied to the actuators and the optimized distances to the collision points are broadcasted to the other agents.}
	\vspace*{-5mm}
\label{fig:archictecture_controlSystem}
\end{figure*}

After addressing the challenge of distributing and efficiently solving the control problem at hand, this section outlines the overall design of the control system architecture, amenable to in-vehicle implementation. A schematic of this architecture, run on a dSPACE MicroAutoBox II, is illustrated in \prettyref{fig:archictecture_controlSystem}. On the input stage, firstly, we receive GNSS positioning and UTC timing information from a low-cost u-blox \mbox{EVK-M8L} GNSS receiver (UBX\_NAV\_PVT message) every \unit[200]{ms}. 
Secondly, we leverage CAN bus signals such as the \chgRevIV{ego vehicle} speed, acceleration and yaw rate, and thirdly, we receive V2V messages from other equipped remote vehicles (RVs). For communication purposes, we utilize a Denso DSRC V2X unit, which is interfaced with the MicroAutoBox through a UDP Ethernet connection. V2V messages involve standardized Cooperative Awareness Messages (CAM) \cite{ETSICAM2014a}, received on average every \unit[100]{ms}, and proprietary Cooperative Control Messages (CCM), which have been designed for the experiment at hand (see \prettyref{sec:archictecture_collPoint}). The CAM messages are stored in a persistent RV buffer to keep track of agents in the vicinity. The input data is then processed by localization algorithms which are in charge of GNSS-based self-localization as well as determining joint collision points with other agents. 
The local MPC controller of Agent $i$ exploits \chgRevIV{the agent's} speed, acceleration and distance to the closest collision point to determine the initial condition $x_k^{[i]}$. Moreover, it utilizes the CCM message, containing the predicted distances $z^{[i]}_{\cdot\mid k}$ of other agents to their collision points, to impose collision avoidance constraints. 
To ensure that the CCM information is valid and consistent when processed within the MPC, every agent's MPC is executed synchronously every \unit[200]{ms}. For that purpose, we make use of a local clock which is synchronized with Universal Time Coordinated (UTC). More precisely, the local clock obtains the UTC time from a GPS message and performs a synchronization step by exploiting an accurate digital UTC trigger signal every UTC second. After optimization, only the first optimal control input $u_{k\mid k}^{[i]\,\star}=a_{x,\text{ref},k\mid k}^{[i]\,\star}$ is applied to the acceleration interface. Moreover, the optimized distances to the joint collision points with other agents are encoded as a CCM message and then broadcasted by the Denso V2X unit. 

To reduce the complexity of localization algorithms, we only contemplate straight crossing agents as indicated in \prettyref{sec:introduction_contribution} and illustrated in \prettyref{fig:archictecture_collPoint_collPointIllustration}. In our experiments, the test driver ensures that the vehicle stays within its designated lane. 
Longitudinal control is taken care of by the control system in \prettyref{fig:archictecture_controlSystem}, which behaves like an \chgAlex{ACC} system that accommodates crossing traffic. While common ACC systems exhibit a sample time less than \unit[100]{ms}, the increased MPC sample time is related to a limited broadcast frequency of the CCM message along with the requirement to solve the underlying OCP in \chgRevIV{real-time}.

\vspace*{-1mm}
\subsection{Self-Localization}
\label{sec:archictecture_selfLocalize}

In our experimental setup, we exploit a low-cost GNSS-based localization system, which receives position updates only with a low frequency (i.e., every \unit[200]{ms}) and asynchronously to the execution of the local MPC controller. To make sure that the latest position information is available whenever the MPC is run and to obtain smooth motion trajectories without discontinuities, we integrate low frequency GNSS measurements with high frequency inertial measurements from the Controller Area Network (CAN) bus by means of an extended Kalman filter (EKF) based estimator, which is run with a sample time of \unit[50]{ms}. In literature, many mature algorithms have already been proposed for this purpose, see \cite{Farrell2008} for a comprehensive overview. These algorithms, though, mostly require \chgRevIV{the ego vehicle's} 
accelerations and angular rates in three dimensions. In our test vehicles, however, we only measure the yaw rate and not the pitch and roll rate. We therefore decided to follow a simpler approach to solve the self-localization problem. Operating in open sky conditions, we are not expecting any GNSS measurement dropouts. That said, our analysis has shown that a simple constant velocity model \cite{Schubert2008a} is sufficient for our use case, i.e., to ensure smooth position trajectories with a high update frequency.
\begin{figure}[b!]
	\centering
	\def\svgwidth{89mm}
	\vspace*{-6mm}
	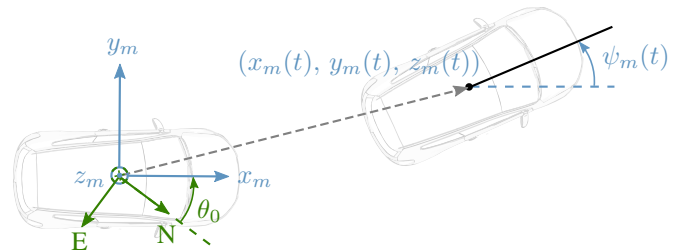
	\vspace*{-7mm}
	\caption{Navigation frame (green) and maneuver reference frame (blue). The origin of the latter refers to the initial pose $P_{\text{WGS},0}$, set at maneuver initiation.}
	\label{fig:archictecture_collisionPoints_coordinateFrames}
\end{figure}

To be more precise, the \textit{Self-Localization} problem at hand aims to estimate the \chgRevIV{ego vehicle's} 
position $(x_m^{[i]}, y_m^{[i]}, z_m^{[i]})$, i.e., the position of its geometric center, and its heading $\psi_m^{[i]}$ in a Cartesian reference frame oriented forward-left-upward which has its fixed origin at the initial maneuver position {--- referred to as \textit{maneuver reference frame}.} 
As we are focusing on the \chgRevIV{ego vehicle}, 
for notational convenience we omit the superscript $[i]$, which has been used in the previous sections to refer to \mbox{Agent $i$.}
Moreover, we subsequently only highlight the most relevant contents of the navigation filter as it is not the main focus of this article.
That said, from the vehicle CAN bus we obtain the vehicle speed $v$ and its yaw rate $\dot\psi$. The GNSS receiver measures the vehicle position 
\(
p_{\text{WGS}} = [\varphi,\, \lambda,\, h]^{\T}
\)
in the World Geodetic System 1984 (WGS-84) in terms of latitude $\varphi$, longitude $\lambda$ and altitude $h$, and the vehicle heading $\theta$ with respect to geographic North, see \cite{Farrell2008} for further details. The underlying constant velocity estimator model approximates the vehicle as a point mass with position $p_m \triangleq [x_m,\, y_m,\, z_m]^{\T}$ and heading $\psi_m$ relative to an initial pose 
\(
P_{\text{WGS},0} \triangleq (\varphi_0,\, \lambda_0,\, h_0,\, \theta_0)
\)
in the WGS-84 coordinate frame, see \prettyref{fig:archictecture_collisionPoints_coordinateFrames} for an illustration. The prediction of states 
\(
x \triangleq [x_m,\, y_m,\, z_m,\, \psi_m ]^{\T}
\)
is performed through numerical integration of 
\begin{align}
\frac{d}{dt}
\underbrace{\begin{bmatrix}
	x_m\\ y_m\\ z_m \\ \psi_m
	\end{bmatrix}}_{x} &= 
\underbrace{\begin{bmatrix}
	v\,\cos(\psi_m)\\ v\,\sin(\psi_m)\\ 0 \\ \dot\psi
	\end{bmatrix}}_{f(x,u)} 
\end{align}
with input vector
\(
u \triangleq [v, \, \dot\psi ]^{\T}.
\) 
We utilize the measurement vector
\(
y \triangleq [\varphi, \, \lambda, \, h, \, \theta]^{\T}
\) 
and derive the measurement equation $y = h_y(x)$ from \cite{Farrell2008}. Essentially, we apply the function 
\begin{align*}
\mathcal{T}_{\text{m2WGS}}: ( x_m,\, y_m,\, z_m,\, \psi_m;\, P_{\text{WGS},0} ) \mapsto ( \varphi,\,  \lambda,\,  h,\,  \theta )
\end{align*}
which transforms the Cartesian estimates $x$ back to the \mbox{WGS-84} frame. 
At the output stage, \textit{Self-Localization} provides the estimated \chgRevIV{ego vehicle} 
pose $P_{m} \triangleq (x_m,\, y_m,\, z_m,\, \psi_m)$ in the maneuver reference frame.

\subsection{Collision Point Estimation}
\label{sec:archictecture_collPoint}
Continuing downstream, \textit{Collision Point Estimation} aims at determining the joint collision point (if it exists) of the \chgRevIV{ego vehicle} 
(i.e., Agent $i$) with other potentially conflicting agents $l \neq i$. As an input, the algorithm consumes the \chgRevIV{ego vehicle} 
pose $P_m^{[i]}$ and kinematic states of other agents $l\neq i$ which are stored in the RV buffer (see \prettyref{fig:archictecture_controlSystem}). RV data contains the WGS-84 position of an \mbox{Agent $l$} along with its heading angle with respect to geographic North, that is, its pose \mbox{$P_{\text{WGS}}^{[l]} \triangleq (\varphi^{[l]}, \lambda^{[l]}, h^{[l]}, \theta^{[l]})$.} 
Through the mapping function
\begin{align*}
\mathcal{T}_{\text{WGS2m}}: ( \varphi,\, \lambda,\, h,\, \theta;\, P_{\text{WGS},0} ) \mapsto ( x_m, \,y_m,\, z_m,\, \psi_m ), 
\end{align*}
we determine Agent $l$'s pose $P_{m}^{[l]}\triangleq(x_{m}^{[l]},\, y_{m}^{[l]},\, z_{m}^{[l]},\, \psi_{m}^{[l]})$ in the \chgRevIV{ego vehicle's} 
maneuver frame.

Before going further into detail, we would like to recall that our particular interest is on scenarios in which the agents are crossing the intersection straight. 
That said, we first predict the future position of Agent $i$ and Agent $l$ at time $t_f > t_k$ \mbox{($t_f$ sufficiently} large), starting at the current time $t_k$, that is, 
\begin{align*}
\hspace{-2mm}
\begin{bmatrix}
x_m(t_f) \\ y_m(t_f)
\end{bmatrix} = 
\begin{bmatrix}
x_m(t_k) \\ y_m(t_k)
\end{bmatrix} + (t_f - t_k)
\begin{bmatrix}
v(t_k) \cos(\psi_m(t_k))\\ v(t_k) \sin(\psi_m(t_k))
\end{bmatrix} 
\end{align*}
with a constant velocity $v$ and constant heading $\psi_m$ and as such with a \chgAlex{yaw rate \mbox{$\dot\psi$} equal to zero} 
due to the agents' straight motion. This way, we obtain the line segment
\begin{align}
\mathcal{S} \triangleq \bigl\{ \,\left( x_m(t_k),\, y_m(t_k)  \right),~ \left( x_m(t_f),\, y_m(t_f)\right)  \,\bigr\}
\label{eq:archictecture_collPoint_lineSegment} 
\end{align}
for Agent $i$ and Agent $l$, which we refer to as $\mathcal{S}^{[i]}$ and $\mathcal{S}^{[l]}$ respectively. Second, we examine whether $\mathcal{S}^{[i]}$ and $\mathcal{S}^{[l]}$ intersect in the horizontal plane (without road inclination, we neglect the $z_m$-dimension), see \prettyref{fig:archictecture_collPoint_collPointIllustration} for an illustration. If there is no intersecting point, there is still the possibility that one of the agents has already passed the collision point. Then, we repeat the calculation by replacing $(x_m(t_k),\, y_m(t_k))$ in \eqref{eq:archictecture_collPoint_lineSegment} with a point from the path history, that is, $(x_m(t_k-t_h),\, y_m(t_k-t_h))$ with sufficiently large $t_h>0$ (first, only for Agent $i$, then only for Agent $l$ and finally for both). Assume, there is an intersecting point, say $( x_{m,\text{CP}}, y_{m,\text{CP}})$. 
Then the distance $\tilde{d}_{c,l}^{[i]} \triangleq s_{c,l}^{[i]}-s^{[i]}$ between Agent $i$'s  joint collision point $s_{c,l}^{[i]}$ with Agent $l$ and its current position $s^{[i]}$ at time $t_k$ is equal to the Euclidean distance between $( x_m(t_k), y_m(t_k)  )$ and $( x_{m,\text{CP}}, y_{m,\text{CP}})$ if Agent $i$ has not yet passed the collision point and equal to the negative Euclidean distance otherwise, see \prettyref{sec:problem_agentDistance}. 
If there is no intersecting point, we define $s_{c,l}^{[i]}=\infty$.

\begin{figure}[t!]
	\centering
	\def\svgwidth{75mm}
	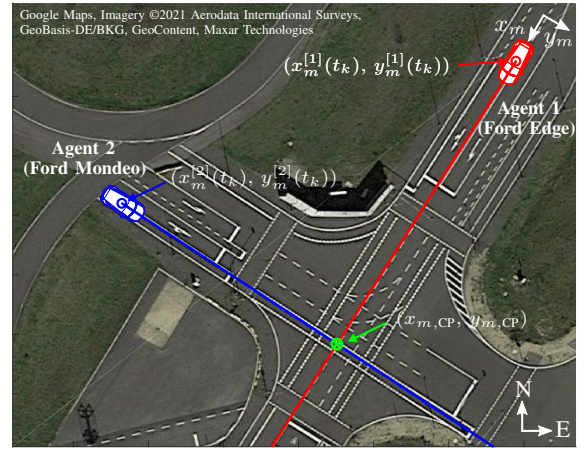
	\vspace*{-2mm}
	\caption{Paths of Agent 1 (red) and Agent 2 (blue), predicted by Agent 1. The straight line predictions intersect in the estimated joint collision point (green).}
	\vspace*{-4mm}
	\label{fig:archictecture_collPoint_collPointIllustration}
\end{figure}

\subsection{Distributed Model Predictive Control}
\label{sec:archictecture_control}

As outlined in the beginning of the section, the local MPC controller is run synchronously on every agent to ensure the synchronicity of broadcasted trajectories $z_{\cdot\mid k}^{[i]}$ and agent decisions $u_{\cdot\mid k}^{[i]}$. For its execution, the MPC requires the initial condition $x_k^{[i]}=[a_{x,k}^{[i]},\, v_{k}^{[i]},\, s_{k}^{[i]}]^{\T}$ and the other agents' distances to the joint collision point with Agent $i$, that is,
\mbox{\(
z_{\cdot \mid k}^{[i]}= \{ (d_{c,i,k+j\mid k}^{[l]} )_{j\in\mathbb{N}_{[1,N]}} \}_{l \in {\mathcal{A}}_{c}^{[i]}}
\)}
as input, see \prettyref{sec:dmpc_solution}.

The initial velocity $v_k^{[i]}$ and acceleration $a_{x,k}^{[i]}$ can directly be obtained from the vehicle CAN bus. The path position $s_{k}^{[i]}$, though, is calculated based on the output of the \textit{Collision Point Estimation} in \prettyref{sec:archictecture_collPoint}. More precisely, it is set to the negative estimated distance to the closest collision point (if it exists), that is, $s_{k}^{[i]} = - \min_{l \in {\mathcal{A}}_{c}^{[i]}}\{ \tilde{d}_{c,l}^{[i]} \}$. If Agent $i$ is not in conflict with any other agent, we set $s_{k}^{[i]}$ to zero. The second input to the MPC, that is, the trajectories $z_{\cdot \mid k}^{[i]}$ can directly be obtained from the CCM message that has been received via V2V communication, see \prettyref{sec:archictecture_CCMmsg}.

After termination of \prettyref{algo:dmpc_solution_CCPappl_CCPalgo}, the optimized distances 
\(
( d_{l,k+1\mid k}^{[i]\,\star}, \ldots, d_{l,k+N\mid k}^{[i]\,\star} )
\)
to every agent $l \in \mathcal{A}_c^{[i]}$
are derived from the solution $u_{\cdot\mid k}^{[i]\,\star}$ of the OCP. When transmitting such distances over V2V to the other agents, though, these distances will be exploited in the next optimization run at time $k+1$. That said, the by one time step shifted trajectory of \mbox{Agent $l$} is required, i.e., beginning at time step $k+2$ until time step $k+N+1$ (when the current time step is $k$). Therefore, we additionally compute ${d}_{l,k+N+1\mid k}^{[i]\,\star}$ by keeping the control input constant after time $k+N-1$, that is, $u_{k+N\mid k}^{[i]} \triangleq u_{k+N-1\mid k}^{[i]\,\star}$. Finally, the trajectories 
\(
\{ d_{l,k+2\mid k}^{[i]\,\star}, \ldots, d_{l,k+N+1\mid k}^{[i]\,\star} \}_{l \in \mathcal{A}_c^{[i]}}
\) 
are forwarded to the \textit{Encode CCM Msg.} block in \prettyref{fig:archictecture_controlSystem}. 

\subsection{Cooperative Control Message (CCM)}
\label{sec:archictecture_CCMmsg}

\chgAlex{With these optimized distances, Agent $i$ sets up the Cooperative Control Message (CCM)}  
\vspace*{-1mm}
\begin{align*}
\text{CCM} \triangleq \left( T_\text{stmp}, \text{ID}^{[i]}, \left\{  \text{ID}^{[l]}, d_{l,k+2\mid k}^{[i]\,\star}, \ldots, d_{l,k+N+1\mid k}^{[i]\,\star} \right\}_{l \in \mathcal{A}_c^{[i]}} \right)\\[-7mm]\notag
\end{align*}
and broadcasts the message to other agents in the vicinity of the intersection.
The CCM contains a time stamp $T_{\text{stmp}}=t_k$ (minutes of hour, milliseconds of minute, 3 Bytes), a unique ID of Agent $i$ (1 Byte), and for every conflicting agent $l$ its respective ID (1 Byte) along with the optimized distances \mbox{(4 Bytes} each, single precision floating point value) of \mbox{Agent $i$} to the joint collision point with Agent $l$ over the prediction horizon of length $N$. 
As such, for every agent $l \in \mathcal{A}_c^{[i]}$, we need to store $1+4N$ Bytes within the message, e.g., 81 Bytes for $N=20$. Additionally, 4 Bytes are required for the time stamp and the \chgRevIV{ego} vehicle ID. 

\begin{figure}[b!]
	\begin{center}
		\vspace*{-2mm}		
		\includegraphics[width=8.0cm]{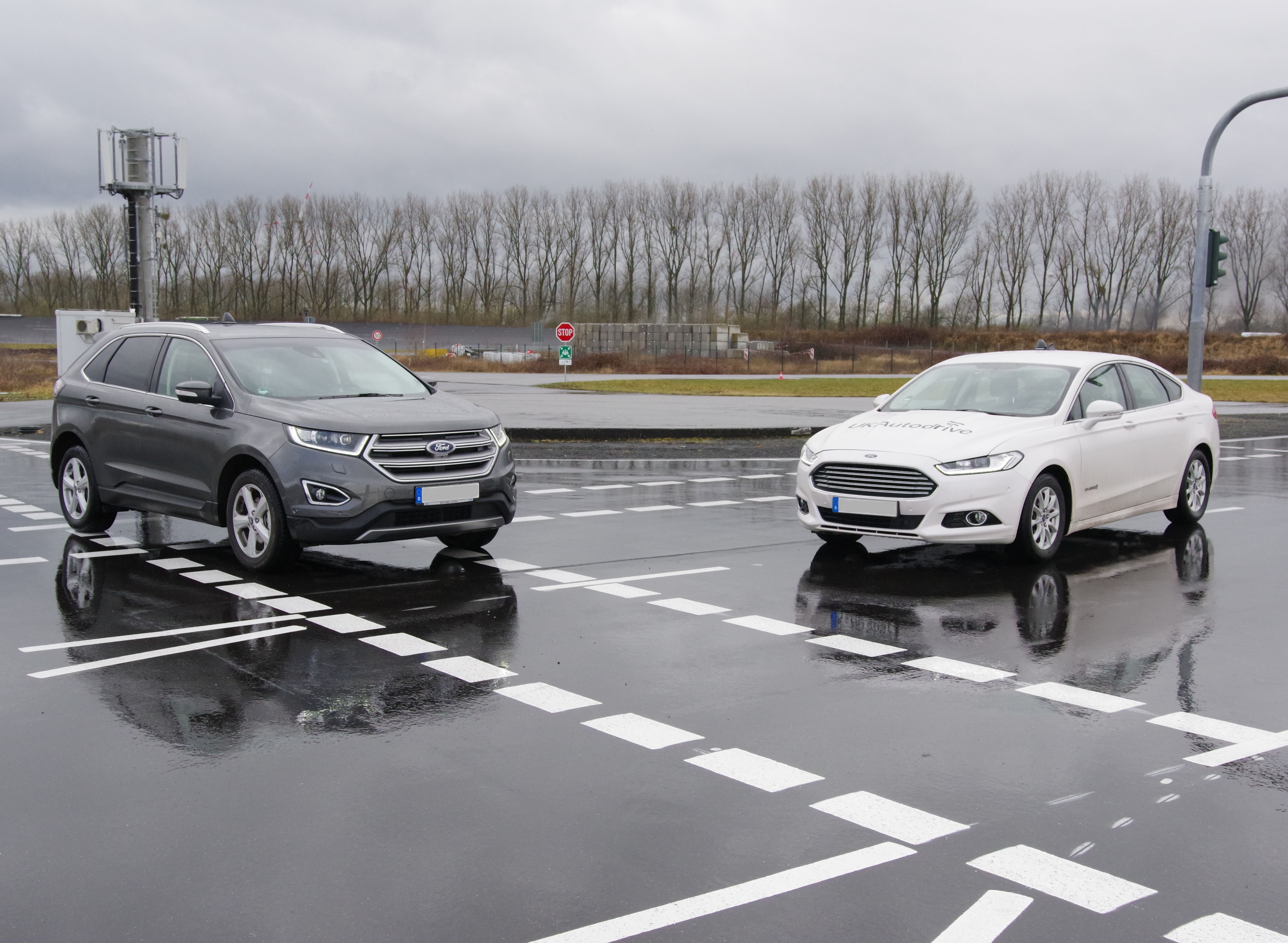}
		\vspace*{-2mm}
		\caption{The two V2X equipped test vehicles on the proving ground:\newline Ford Edge (left) is Agent 1 and Ford Mondeo Hybrid (right) is Agent 2.} 
		\vspace*{-4mm}
		\label{fig:expResults_setup_testvehicles}
	\end{center}
\end{figure}

\section{Experimental Results}
\label{sec:expResults}

\subsection{Experimental Setup and Parameterization}
\label{sec:expResults_setup}

To evaluate the proposed control system architecture in experimental tests, two test vehicles have been available to the research team, that is, a Ford Edge (\mbox{Agent 1}) and a Ford Mondeo Hybrid (\mbox{Agent 2}), see \prettyref{fig:expResults_setup_testvehicles}. Each test vehicle is equipped with the same hardware and software as described in \prettyref{sec:archictecture}. A centimeter-precision RTK positioning system has not been available for our experiments such as to serve as ground truth with respect to the low-cost GNSS used for control purposes. The main aim of the control system, though, is to satisfy collision avoidance constraints given the low-cost GNSS positioning information. So, the availability of ground truth measurements would help to assess the actual distance between agents when crossing the intersection but is not necessarily required in our case. Actually, we never encountered critical situations due to position inaccuracies, in which the agents got too close. To accommodate positioning uncertainties within the control system, we have anyway added an additional error budget to the minimum distance between the agents. 
The test drives have been carried out on the Aldenhoven Testing Center close to Aachen, Germany. We utilized \chgAlex{a single lane four-way intersection which has been crossed straight by the agents.} With the given control system architecture and the experimental setup, all fundamental assumptions that have been made in \prettyref{sec:problem_description} are satisfied.

\begin{table}[b!]
	\vspace*{-4mm}
	\begin{center}
		\caption{Parameterization of the local MPC controllers\vspace*{-6mm}}
		\label{tab:expResults_setup_MPCparams}
		\begin{tabular}{lccc}
			\hline		
			& \hspace*{-1mm}\textbf{Agent 1}	& \hspace*{-1.5mm}\textbf{Agent 2} 	\\ 
			& \hspace*{-1mm}(Ford Edge)	& \hspace*{-1.5mm}(Ford Mondeo) 	\\ \hline\\
			\vspace*{-0.5cm}&&\\
			\multicolumn{3}{l}{\textit{\textbf{Common Parameters}}}\\			
			MPC Sample Time $[\unit[]{s}]$ ~~~ & 0.2 & 0.2 \\			
			Horizon Length $N$ $[-]$ ~~~ & 20 & 20 \\			
			Weights $(Q,\,Q_N,\,R,\,S)$ ~~~  & (1, 1, 5, 5) & (1, 1, 5, 5)\\		
			Min./Max. Ref. Accel. $[\unit[]{m/s^2}]$ ~~~ & (-5, 2) & (-5, 2) \\								
			Safety Distance $d_{\text{safe}}$ $[\unit[]{m}]$ ~~~  & 15 & 15 \\	
			Vehicle Length $L$, Width $W$ $[\unit[]{m}]$ ~~~  & (4.8, 1.9) & (4.8, 1.9) \\					
			Priority $[\unit[]{-}]$                                  ~~~ &  2 &  1 \\		
			
			\vspace*{-0.2cm}&&\\
			\multicolumn{3}{l}{\textit{\textbf{Scenario 1:} Lower Urban Speed Limit}}\\
			Ref. Speed $v_{\text{ref}}$, Max. Speed
		    $[\unit[]{m/s}]$ ~ & (12, 13.2) & (10, 11)\\			
			Init. Condition\chgRevI{\textsuperscript{1}} $(s_0\,[\unit[]{m}],\,v_0\,[\unit[]{m/s}])$	  ~~~ &(-83.5, 11.9)  & (-64.8, 10.0) \\
			\vspace*{-0.15cm}&&\\
			\multicolumn{3}{l}{\textit{\textbf{Scenario 2:} Regular Urban Speed Limit}}\\
			Ref. Speed $v_{\text{ref}}$, Max. Speed
$[\unit[]{m/s}]$ ~ & (15, 16.5) & (11, 12.1)\\			
			Init. Condition\chgRevI{\textsuperscript{1}} $(s_0\,[\unit[]{m}],\,v_0\,[\unit[]{m/s}])$	  ~~~  & (-103.1, 14.8) & (-66.7, 10.3) \\
			\hline\\[-2mm]
			\multicolumn{3}{l}{\chgRevI{\textsuperscript{1} \scriptsize Initial distance to collision point is equal to $|s_0|$.}}
		\end{tabular}
	\end{center}
	\vspace*{-5mm}
\end{table}
Subsequently, we discuss two scenarios that mostly differ with respect to the agents' reference speeds. Particularly, in \mbox{Scenario 1} the agents exhibit a set-speed of $\unit[12]{m/s}$ \mbox{(Agent 1)} and $\unit[10]{m/s}$ (Agent 2) respectively. These speeds can, e.g., be observed in urban areas with lower urban speed limit of $\unit[30]{kph}$. Then, in \mbox{Scenario 2} we increase the reference speed of Agent 1 to $\unit[15]{m/s}$ which is slightly above the regular urban speed limit of $\unit[50]{kph}$ in Germany. Due to a limited track length, though, the reference speed of \mbox{Agent 2} is almost the same as in Scenario 1. 
For both scenarios, the agents' maximum speed is set to $\overline{v}^{[i]}=1.1 v_{\text{ref}}^{[i]}$} and their initial configuration corresponds to \prettyref{fig:archictecture_collPoint_collPointIllustration}. 
Moreover, \mbox{Agent 2} always exhibits a higher priority than Agent 1. Recall that the definition of priority defines which agent needs to impose collision avoidance constraints rather than an intersection crossing order.
Nonetheless, we have intentionally chosen the initial conditions such that \mbox{Agent 1} always has to brake for Agent 2. Otherwise, both agents would exhibit almost no control action, making it less attractive for our analysis. 

In terms of parameterization, Agent 1 and Agent 2 have the same length and width of $L=\unit[4.8]{m}$ and $W=\unit[1.9]{m}$ respectively. To ensure safe intersection crossing, the minimum safety distance $d_{\text{safe}}$, defined as the distance between the agents' geometric centers along their path coordinate $s$ (see \prettyref{sec:problem_agentDistance}), has been set to $\unit[15]{m}$. In the perpendicular straight crossing scenario, the minimum distance between the vehicle bounding boxes amounts to $d_{\text{safe}} - (W+L)/2 = \unit[11.65]{m}$. 
For the dynamic powertrain time constant $T_{a_x}^{[i]}$, we recognized during system identification that this quantity \chgRevI{evolves as a function} of the current state $x^{[i]}$ and input $u^{[i]}$. To this end, the MPC utilizes a lookup table $T_{a_x}^{[i]}(x_k^{[i]},u_k^{[i]})$ to determine the time constant at every time step $k$.
The local MPC controllers exhibit a sample time of $\unit[200]{ms}$ with a prediction horizon of 20 steps, thus covering a preview time of $\unit[4]{s}$. To solve OCP \eqref{eq:dmpc_solution_CCPappl_CCPalgo_convexOCP} on the dSPACE \mbox{MicroAutoBox II}, qpOASES \cite{Ferreau2014} is applied as QP solver. If terminal constraint \eqref{eq:cmpc_objCons_vMinMeanCons} renders OCP \eqref{eq:dmpc_solution_CCPappl_CCPalgo_convexOCP} infeasible, we perform an appropriate braking maneuver to prevent the agent from entering the intersection. \prettyref{tab:expResults_setup_MPCparams} conveys the most relevant parameters.

\begin{figure}[b!] 
	\setlength\fwidth{0.42\textwidth}		
	\vspace*{-6mm}
	\input{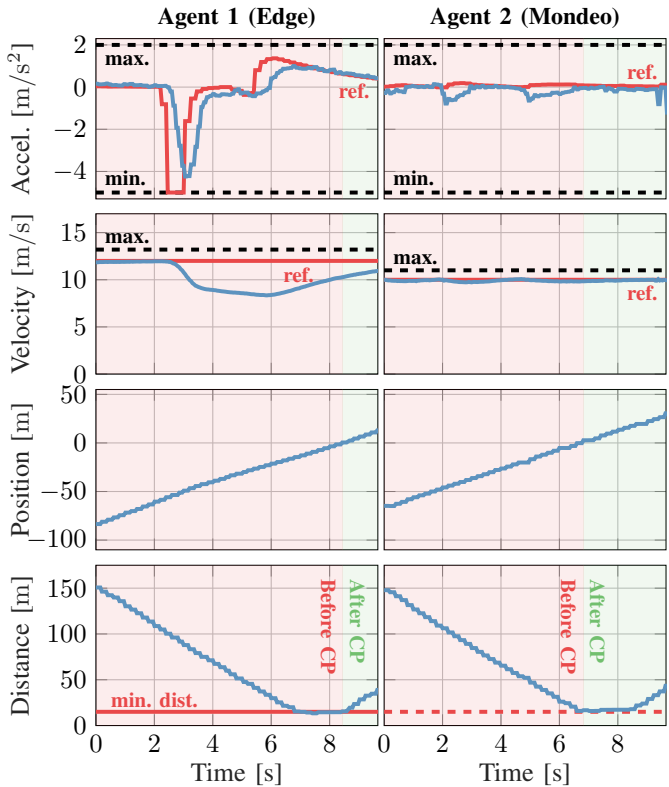}
	\vspace*{-8mm}
	\caption{Scenario 1 (Lower Urban Speed Limit): Acceleration, velocity and path position of both agents along with their distance to each other. Agent 1 (low priority) is able to satisfy collision avoidance constraints at all times.}
	\label{fig:results_discussion_scenario1Timeplot}
\end{figure}

\subsection{Discussion of Results}
\label{sec:expResults_discussion}

\subsubsection{Scenario 1 (Lower Urban Speed Limit)}

\begin{figure}[b!] 
	\setlength\fwidth{0.42\textwidth}		
	\vspace*{-6mm}	
	\input{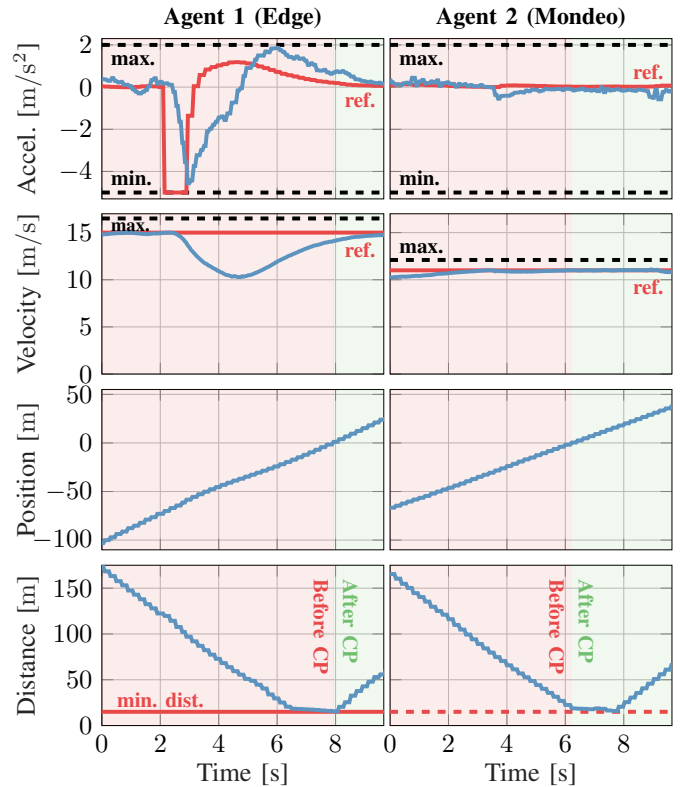}
	\vspace*{-8mm}
	\caption{Scenario 2 (Regular Urban Speed Limit): Acceleration, velocity and path position of both agents along with their distance to each other. Agent 1 (low priority) is able to satisfy collision avoidance constraints at all times.}
	\label{fig:results_discussion_scenario2Timeplot}
\end{figure}

\prettyref{fig:results_discussion_scenario1Timeplot} illustrates the experimental results for Scenario 1. With  Agent 1 in the left and Agent 2 in the right column, the figure reveals from top to bottom: 
1) actual (solid blue) and reference acceleration (solid red) along with the upper and lower bound of the reference (dashed black); 2) actual (solid blue), reference (solid red) and maximum velocity (dashed black); 3) actual path position (solid blue). The bottom plot provides the distance \chgRevI{$\mathrm{dist}(1,2)$ and $\mathrm{dist}(2,1)$ (see \eqref{eq:problem_agentDistance_distAgentToAgent})} between the agents \chgRevI{along their paths}, derived from the calculations in \prettyref{sec:archictecture_collPoint}, along with the minimum safety distance that has to be ensured by \mbox{Agent 1}. As Agent 2 exhibits higher priority, the dashed red line only indicates the minimum safety distance while collision avoidance constraints are not imposed on this agent. 
\chgRevI{It should be noted that $\mathrm{dist}(1,2)$ and $\mathrm{dist}(2,1)$ are the same in theory. In practice, though, they 
depend on every agent's local GNSS measurements, the estimated collision point and the optimized distance of the other agent to the respective collision point, see \prettyref{sec:archictecture_control}. For that reason, we have plotted both quantities to provide evidence that a safe distance is indeed ensured for both agents.}
Moreover, the red patches indicate the time interval when the respective agent is approaching the joint collision point (CP), i.e. $s^{[i]} \leq s_{c,l}^{[i]}$, while the green patches highlight the time interval when the agent is moving away from it, i.e., $s^{[i]} > s_{c,l}^{[i]}$.

To start with, in the first scenario Agent 1 (low priority) exhibits an initial speed of $\unit[11.9]{m/s}$ while Agent 2 (high priority) approaches the collision point with $\unit[10]{m/s}$ --- corresponding to speeds that can be observed in lower urban speed limit areas. 
By evidence of \prettyref{fig:results_discussion_scenario1Timeplot}, Agent 2 crosses the intersection with constant velocity and without any reaction to Agent 1. Minor acceleration demands results from a slight inclination of the road section. Conversely, Agent 1 starts to decelerate at $t=\unit[2.1]{s}$ to give right of way to Agent 2. More precisely, Agent 1 decelerates with a maximum deceleration of $\unit[-4.2]{m/s^2}$. To compensate the drivetrain lag, even higher decelerations of up to $\unit[-5]{m/s^2}$ are requested by the MPC. This way, Agent 1 slows down to a minimum speed of $\unit[8.3]{m/s}$. After Agent 2 has passed the joint collision point at $t=\unit[6.8]{s}$, Agent 1 can safely pass that point at $t=\unit[8.4]{s}$ and continues to track its reference speed. Most importantly, during the entire maneuver a minimum safety distance of $\unit[15]{m}$ can be guaranteed by Agent 1 (with respect to the low-cost GNSS positioning information). That said, the proposed control concept successfully accomplishes the maneuver. 

\begin{figure}[t!] 
	\setlength\fwidth{0.44\textwidth}		
	\input{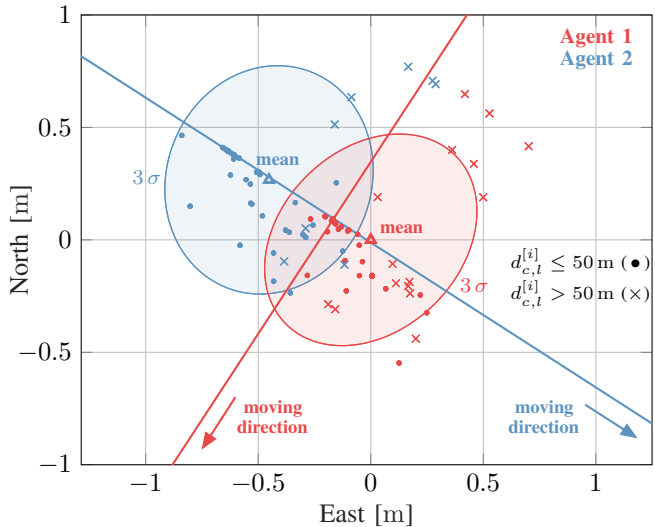}
	\vspace*{-4mm}
	\caption{Estimated collision points of both agents for Scenario 1. 
	At distances $\leq\unit[50]{m}$, the collision points are shown as a filled circles ($\bullet$), and as crosses ($\times$) at larger distances. Additionally, we illustrate the mean (triangle) and three times the standard deviation (circle). The solid lines are the \mbox{agents' paths.}\vspace*{-4mm}} 
	\label{fig:results_discussion_scenario1cpoint}
\end{figure}

Besides the satisfaction of control objectives and constraints, we have additionally analyzed the variation of the estimated collision point, required to determine whether Agent $i$ is in conflict with Agent $l$ as well as to compute its initial condition $x_k^{[i]}$ at time $t_k$.
Especially at larger distances, the estimated position of the collision point is very sensitive to small changes of the agent's heading angle and as such to GNSS heading errors or driver steering inputs. \prettyref{fig:results_discussion_scenario1cpoint} provides an overview of the distribution of estimated collision points in a Cartesian (North, East) coordinate frame for Scenario 1. We represent all collision points at a distance less than or equal to \unit[50]{m} with a filled circle ($\bullet$) and above that threshold with a cross ($\times$). Moreover, assuming a normal distribution $\mathcal{N}(\mu,\sigma)$ we show the mean $\mu$ ($\triangle$) along with the $3\sigma$ standard deviation. Finally, the solid lines highlight the agents' paths. It can be recognized that the maximum deviation of the estimated collision points from the actual one is less than $\unit[1]{m}$ during the entire maneuver and even less when the agents get closer to that point. Our experiments have shown that our control system is robust to variations of that magnitude. If GNSS errors increase or the driver steering input is inappropriate (an uncontrollable noise factor to our system), larger variations may cause an uncomfortable driving behavior. When, e.g., the distance to the collision point suddenly gets smaller, the control system may need to decelerate more severely to satisfy collision avoidance constraints. To make the control system robust against such noise factors has actually not been in the scope of this work.

\subsubsection{Scenario 2 (Regular Urban Speed Limit)}
The second scenario mainly differs with respect to the agents' speed. That said, Agent 1 (low priority) exhibits an initial speed of $\unit[14.8]{m/s}$ which corresponds to a regular urban speed limit. Agent 2 (high priority) approaches the collision point with an initial speed of $\unit[10.3]{m/s}$. 
Similar to Scenario 1, Agent 2 passes the intersection with a constant velocity and without the need to react to Agent 1. With an increased speed, \mbox{Agent 1} needs to decelerate more heavily compared to Scenario 1, that is, at $t=\unit[2.0]{s}$ with a maximum deceleration of $\unit[-4.6]{m/s^2}$. As a consequence, Agent 1 slows down to $\unit[10.3]{m/s}$ to let \mbox{Agent 2} pass the joint collision point at $t=\unit[6.2]{s}$. After \mbox{Agent 1} has crossed the intersection, Agent 2 follows at $t=\unit[7.9]{s}$ and resumes to track its reference speed. In spite of the higher maneuver speed, it is evident that collision avoidance constraints can still be satisfied at all times. 

To conclude, the proposed control system architecture has successfully been evaluated in urban driving scenarios. More precisely, it accommodates control objectives, ensures collision avoidance and is amenable to in-vehicle implementation. The latter statement is further supported by the fact that all calculations have been finished within the given sample time on the dSPACE MicroAutoBox II.

\section{Conclusion}
\label{sec:conclusion}
We have conveyed a fully distributed control system architecture to safely coordinate CAVs \chgRevIV{at} road intersections with no traffic signs or lights. For control purposes, a fully distributed MPC scheme has been proposed. To allow every agent to solve its originally nonconvex OCP fast, penalty CCP is applied to obtain a local solution in \chgRevIV{real-time}. 
For an in-vehicle implementation, the control layer is complemented with a localization layer to estimate the agents' positions and their joint collision points. The entire control system architecture is implemented on two test vehicles and the respective algorithms are run on a dSPACE MicroAutoBox II. Two experimental tests, i.e., a lower urban speed limit and a regular urban speed limit scenario have demonstrated that the proposed concept satisfies control objectives and ensures collision avoidance. 

\chgAlex{While agent priorities have been fixed in this work, our current research works \cite{Molinari2020a} also investigate time-varying priorities. 
Moreover, we aim to solve the centralized OCP through distributed numerical optimization methods w/o prioritization.}


%



%
%

\ifCLASSOPTIONcaptionsoff
  \newpage
\fi



\vspace*{-1mm}
\bibliographystyle{IEEEtran}
\bibliography{IEEEabrv}
%
%
%

%

\vspace*{-10mm}
\begin{IEEEbiography}[{\includegraphics[width=1in,height=1.25in,clip,keepaspectratio]{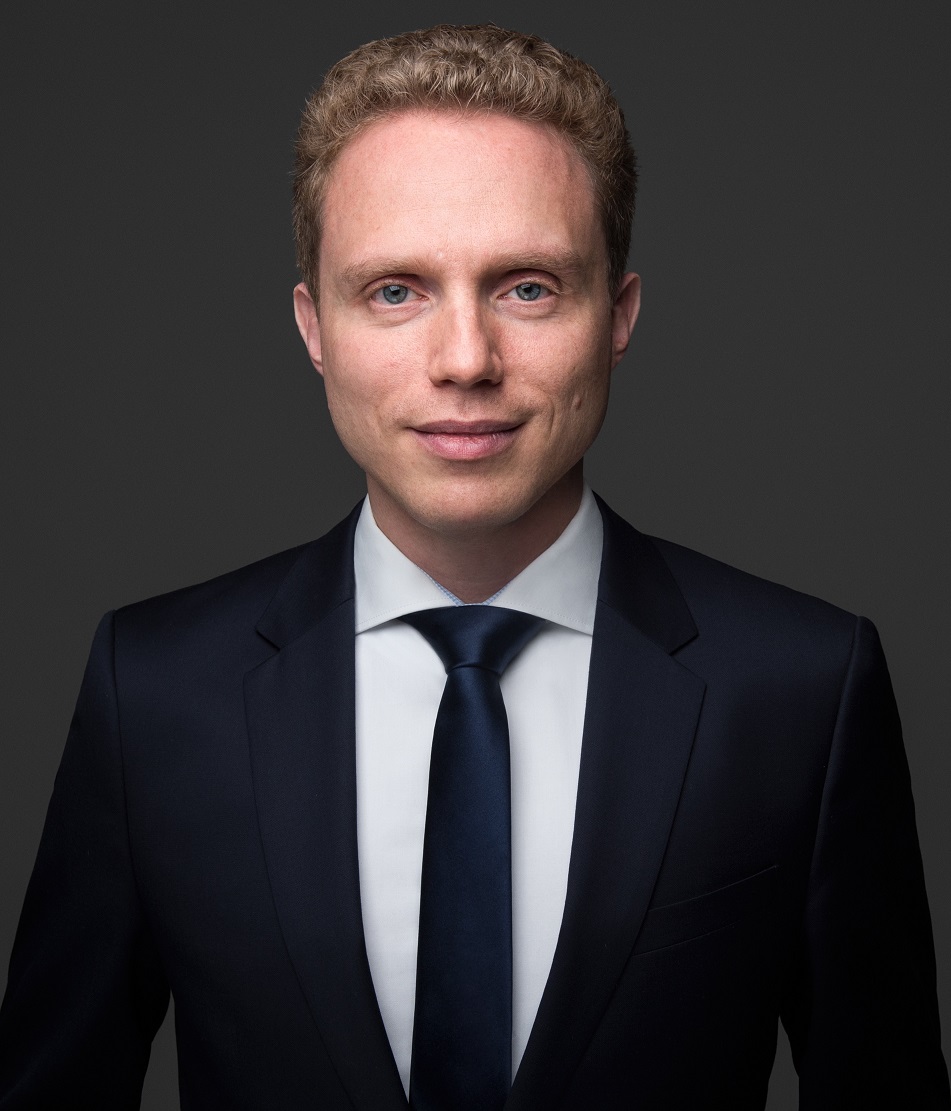}}]{Alexander Katriniok} (M'15--SM'19) received the PhD in Mechanical Engineering from RWTH Aachen University, Aachen, Germany, in 2013. 	Since 2016, Dr. Katriniok is with the Ford Research \& Innovation Center (RIC) in Aachen, Germany. He is working on sensing \& perception, machine learning and advanced control methods for connected and automated driving applications. His scientific research interests include learning-/data-based MPC, distributed optimal control and (distributed) numerical optimization with application to motion planning and control of automated vehicles and robots in uncertain environments. 	
\end{IEEEbiography}
\vspace*{-10mm}
\begin{IEEEbiography}[{\includegraphics[width=1in,height=1.25in,clip,keepaspectratio]{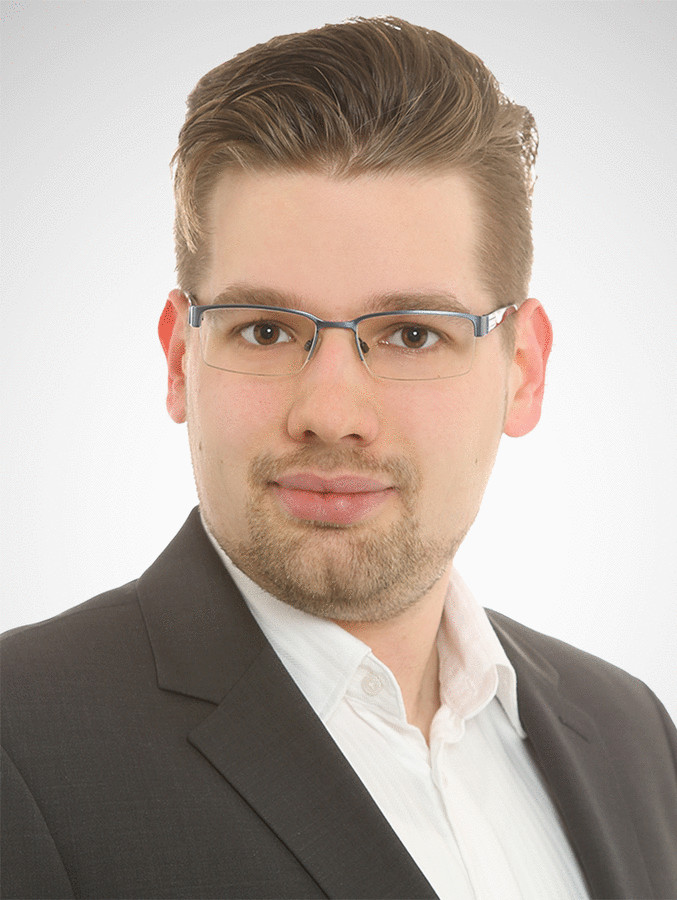}}]{Benedikt Rosarius}
received the master degree in Electrical Engineering from RWTH Aachen University, Aachen, Germany, in 2019. 
As part of his master thesis at Ford, he focused on V2V-based intersection automation utilizing distributed model predictive control. His research focused on extending and applying the distributed MPC algorithm to in-vehicle experiments, and the comparison of simulation-based and experimental results in on-road scenarios.
\end{IEEEbiography}
\vspace*{-10mm}
\begin{IEEEbiography}[{\includegraphics[width=1in,height=1.25in,clip,keepaspectratio]{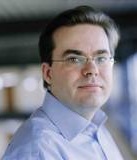}}]{Petri~M\"{a}h\"{o}nen} (SM'01)
is currently a Full Professor and the Chair of Networked Systems with RWTH Aachen
University. His current research focuses on cognitive radio systems, embedded intelligence, future wireless networks architectures, including MillimeterWave systems, and techno-economics especially from a regulatory perspective. He is currently serving as Editor of IEEE Transactions of Wireless Communications.
\end{IEEEbiography}





\end{document}